\documentclass[reqno,oneside]{amsart}
\usepackage{amsfonts,cite,amssymb,geometry,amsmath,amsthm}
\usepackage[all]{xy}

\newtheorem{thm}{Theorem}[section]
\newtheorem{cor}[thm]{Corollary}
\newtheorem{lem}[thm]{Lemma}
\newtheorem{prop}[thm]{Proposition}
\theoremstyle{remark}
\newtheorem{rem}[thm]{Remark}
\numberwithin{equation}{section}

\def\sz{\sigma}

\begin{document}

\title[] {Average and probabilistic linear widths of sequence spaces  with a Gaussian measure}

\author[]{Yun Ling} \address{School of Mathematics and Statistics, Hanjiang Normal University, Shiyan 442000, Hubei, China}
\email{lingyun@hjnu.edu.cn}

\author[]{Heping Wang} \address{School of Mathematical Sciences, Capital Normal
University, Beijing 100048,
China}
\email{wanghp@cnu.edu.cn}


\begin{abstract}
Let $\gamma$ be a Gaussian measure on the sequence space $\ell_2$
whose mean is zero and whose correlation operator $C_\gamma$ is a
diagonal operator on $\ell_2$. Exact values of the $q$-average
linear $n$-widths in $\ell_q$ for $1\le q<\infty$ are obtained. We
also obtain   asymptotic orders of  the average linear $n$-width
in $\ell_\infty$ under some condition on $C_\gamma$,  and of the
probabilistic linear $(n,\delta)$-widths in $\ell_q$ for $1\le
q\le 2$, and for $2<q\le \infty$ under some condition on
$C_\gamma$.
\end{abstract}
\maketitle

\section{Introduction}
Let $K$ be a bounded subset of a normed linear space $X$ with norm $\|\cdot\|_{X}$. The classical linear and Kolmogorov $n$-widths of the set $K$ in $X$ are defined by
$$\lambda_{n}(K,X):=\inf_{L_{n}}\sup_{x\in K}\|x-L_{n}x\|_{X},$$
and
$$d_{n}(K,X):=\inf_{F_{n}}\sup_{x\in K}\inf_{y\in F_{n}}\|x-y\|_{X},$$
respectively, where $L_{n}$ runs over all linear operators from
$X$ to $X$ with rank at most $n$, and $F_{n}$ runs through all
possible linear subspaces of $X$ with dimension at most $n$.
Detailed information about the classical widths, such as the linear
and Kolmogorov widths, may be found in \cite{H, LGM, PI, P}.

Next we define approximation in the average and probabilistic case setting. Let
$W$ be a separable Banach space and assume that $W$ contains a
Borel field $\mathcal{B}$ consisting of open subsets of $W$ and is
equipped with a probability measure $\mu$ defined on
$\mathcal{B}$. For $0<p<\infty$, the $p$-average linear and
Kolmogorov $n$-widths are defined by
$$\lambda^{(a)}_{n}(W,\mu,X)_{p}:=\inf_{L_{n}}
\Big(\int_{W}\|x-L_{n}x\|^{p}_{X}d\mu(x)\Big)^{1/p},$$ and
$$d^{(a)}_{n}(W,\mu,X)_{p}:=\inf_{F_{n}}
\Big(\int_{W}\inf_{y\in F_{n}}\|x-y\|^{p}_{X}d\mu(x)\Big)^{1/p},$$
respectively. They reflect the optimal approximations of ``most''
elements of classes by $n$-dimensional subspaces and linear
operators with rank $n$. Therefore, in the worst case setting, the
approximation emphasizes the behavior of the ``worst'' elements,
while in the average case setting, the emphasis is on
the behavior of ``most'' elements.

Let $\delta\in[0,1)$. The probabilistic linear and Kolmogorov
$(n,\delta)$-widths of a set $W$ with a measure $\mu$ in the space
$X$ are defined by
$$\lambda_{n,\delta}(W,\mu,X):=\inf_{G_{\delta}}\lambda_{n}(W\backslash G_{\delta},X),$$
and
$$d_{n,\delta}(W,\mu,X):=\inf_{G_{\delta}}d_{n}(W\backslash G_{\delta},X),$$
respectively, where $G_{\delta}$ runs through all possible
measurable subsets of $W$  with measure
$\mu(G_{\delta})\leq\delta$. Hence, the probabilistic linear and
Kolmogorov $(n,\delta)$-widths can be understood as the
$\mu$-distribution of the approximation on all measurable subsets
of $W$ by linear operators with rank $n$ and by $n$-dimensional
subspaces. Therefore, the probabilistic case setting reflects the
intrinsic structure of the class, and compared with the worst case
setting, allows one to give a deeper analysis of the approximation
for the function class.

There are many papers devoted to investigating the average and
probabilistic linear and Kolmogorov widths, and a great number of
interesting results have been obtained  (see \cite{CNL, DW, FY,
M1, M2, M3, R, TWW, W1, W2}). In this paper, we discuss the
average and probabilistic linear widths of sequence spaces with a
Gaussian measure. For $1\le N\le \infty$, let $\ell_{q}^N,\,1\leq
q\leq\infty$, be the classical real sequence space with the usual
norm
$$\|x\|_{q}\equiv\|x\|_{\ell_{q}^N}:=\left\{
\begin{aligned}
 &\Bigg(\sum^{N}_{j=1}|x_{j}|^{q}\Bigg)^{1/q},&&{\rm if}\ 1\leq q<\infty,\\
 &\max_{1\le j\leq N}|x_{j}|,&&{\rm if}\ q=\infty,
\end{aligned}
\right.$$ and let $B_{q}^N:=\{x\in\ell_{q}^N:\|x\|_{q}\leq1\}$
denote the unit ball of $\ell_{q}^N$. As usual, $\ell_{2}^N$ is a
Hilbert space with respect to the inner product $\langle
x,y\rangle:=\sum^{N}_{j=1}x_{j}y_{j}$. Let $\{e_j\}_{j=1}^N$ be
the canonical orthonormal basis in $\ell_2^N$, where
$e_j=\{\delta_{jk}\}_{k=1}^N$,  and $\delta_{jk}$ is  the
Kronecker delta
 $$\delta_{jk}:=\left\{
\begin{aligned}
 &1,&&j=k,\\
 &0,&&j\neq k.
\end{aligned}
\right.$$ For simplicity, we write $\ell_q$ and $B_q$ instead
of $\ell_q^\infty$ and $B_q^\infty$, respectively.

We equip $\ell_2^N$ with a Gaussian measure $\gamma$ whose mean is
zero and whose correlation operator $C_\gamma$ has eigenvectors
$e_j,\  j=1,\dots, N$ and eigenvalues $\sigma_j, \ j=1,\dots, N$.
We assume that $$\sigma_1\ge \sigma_2\ge \dots \ge \sigma_N\ge 0,
\ \ {\rm and}\ \ \sum_{j=1}^N\sigma_j<+\infty.$$ Clearly, the
correlation operator $C_\gamma$ is the diagonal operator ${\rm
diag}(\sigma_1, \dots, \sigma_N)$, i.e.,
$$C_\gamma e_j=\sigma_j e_j,\ \ j=1,\dots, N,$$
and
$$C_\gamma x=\sum_{j=1}^N\sigma_j x_j e_j \ \ {\rm for}\ \ x\in\ell_2^N. $$
We remark that $\gamma$ is a Gaussian measure on $\ell_2$ if and
only if the trace of the correlation operator $C_\gamma$ is finite
(see \cite[Theorem 2.3.1]{Bo}), i.e.,
$$ {\rm trace}\ C_\gamma= \sum_{j=1}^\infty \sigma_j<+\infty.$$
For $1\le n<N$, denote by  $S_n$
the Fourier partial summation operator on $\ell_2^N$ defined by
$$S_n x:=\sum_{j=1}^nx_je_j=(x_1,\dots, x_n,0,\dots).$$ If
$N=m<+\infty$ and $\sigma_1= \dots=\sigma_m=1$, then $\gamma$ is
just the standard Gaussian measure $\gamma_m$ on $\mathbb{R}^{m}$, given by
$$\gamma_{m}(G)=(2\pi)^{-m/2}\int_{G}\exp(-\|x\|^{2}_{2}/2)dx,\quad {\rm for\ each\ Borel\ subset} \ G\subset\mathbb{R}^{m}.$$

The theory of widths has important applications in various fields such as harmonic analysis, functional analysis and operator theory, and numerical analysis and information-based complexity (see \cite{Be,DD,FY,FY2004,Hi,La,MR,M1,P1990,P1993,S1959,S1960,SS,Vo,Wo,Y}).

This paper is devoted to investigating average and probabilistic
linear widths in $\ell_{q}$ with the Gaussian measure $\gamma$.
There is an extensive literature on average and probabilistic
linear widths of $\mathbb{R}^{m}$ with the standard Gaussian
measure $\gamma_m$ in $\ell_{q}^{m}\ (1\le q\le \infty)$. For
example, the authors in \cite{FY} have obtained  exact orders of
the probabilistic linear $(n,\delta)$-widths
$\lambda_{n,\delta}(\mathbb{R}^{m},\gamma_{m},\ell^{m}_{q})$ for
$1\leq q<2,\ m\geq2n$, and the $p$-average linear $n$-widths
$\lambda^{(a)}_{n}(\mathbb{R}^{m},\gamma_{m},\ell^{m}_{q})_{p}$
for $1\leq q<\infty,\ 0<p<\infty,\ m\geq2n$. The authors in
\cite{M2,M3} have obtained  exact orders of
$\lambda_{n,\delta}(\mathbb{R}^{m},\gamma_{m},\ell^{m}_{q})$ for
$2\leq q\leq\infty,\ m\geq2n$, and
$\lambda^{(a)}_{n}(\mathbb{R}^{m},\gamma_{m},\ell^{m}_{q})_{p}$
for $q=\infty,\ 0<p<\infty,\ m\geq2n$. The authors in \cite{CNL}
have obtained  exact orders of
$\lambda_{n,\delta}(M:\mathbb{R}^{m}\rightarrow\ell^{m}_{q},\gamma_{m})$
and
$\lambda^{(a)}_{n}(M:\mathbb{R}^{m}\rightarrow\ell^{m}_{q},\gamma_{m})_{p}$
for $ 2\leq q<\infty,\ 0<p<\infty,\ m\geq2n$ of a certain diagonal
matrix $M$. Later, the authors in \cite{DW} obtained exact values
of
$\lambda^{(a)}_{n}(D:\mathbb{R}^m\rightarrow\ell_{q}^m,\gamma_{m})_{q}$,
$1\le q<\infty$, and asymptotic estimates of
$\lambda_{n,\delta}(D:\mathbb{R}^m\rightarrow\ell_{q}^m,\gamma_m),\,
1\le q\le2$ for diagonal matrices $D$. We obtain exact values of
the $q$-average linear $n$-widths in $\ell_q$ for $1\le q<\infty$.
Additionally, we determine  asymptotic orders of  the average
linear $n$-width in $\ell_\infty$ under some condition on
$C_\gamma$,  and of the probabilistic linear $(n,\delta)$-widths
in $\ell_q$ for $1\le q\le 2$, and for $2<q\le \infty$ under some
condition on $C_\gamma$. Our contributions can be summarized in
three aspects.

Firstly, we extend  results on average and probabilistic widths
with respect to Gaussian measures from finite-dimensional spaces
to infinite-dimensional spaces. The extension requires taking into
account a fundamental difference between the Gaussian measure on
$\mathbb{R}^{m}$ and that on $\ell_{2}$. In $\mathbb{R}^{m}$,
 the standard Gaussian measure $\gamma_m$ can be represented by a density
function with respect to  the translation-invariant Lebesgue
measure. Consequently, the covariance operator of $\gamma_m$ can
be the identity operator on $\mathbb R^m$. However, in $\ell_{2}$
no such nontrivial, $\sigma$-finite translation-invariant measure
exists. Consequently, a Gaussian measure  cannot be represented by
a density function with respect to any translation-invariant
measure. Moreover, the covariance operator of any Gaussian measure
on $\ell_{2}$ cannot be the identity operator. Hence, our proof
requires a more refined method.

Secondly, by combining the comparison theorem for Gaussian
measures, classical linear widths, and the result of Maiorov and
Wasilkowski in \cite{M3}, we obtain  exact orders of the average
linear widths under the $\ell_{\infty}$-norm when the correlation
operator $C_{\gamma}$ satisfies certain conditions.

Thirdly, by employing classical and average linear widths,
together with tail estimates of the standard Gaussian distribution
and concentration inequalities for centered Gaussian random
vectors on $\ell_{2}$, we obtain  asymptotic orders of the
probabilistic linear $(n,\delta)$-widths in $\ell_{q}\
(2<q\leq\infty)$ under certain conditions on $C_\gamma$.

Our first theorem provides  exact values of
average linear widths
$\lambda^{(a)}_{n}(\ell_{2},\gamma,\ell_{q})_{q}$ for $1\leq
q<\infty$.

\begin{thm}\label{thm1.1}
Let $1\le q<\infty$ and $C(q)=\pi^{-\frac{1}{2q}}
2^{\frac{1}{2}}\,\Gamma\left(\frac{q+1}{2}\right)^{\frac{1}{q}}$.

{\rm(i)} If $1\leq q\leq 2$ and
$\sum^{\infty}_{j=1}\sigma^{q/2}_{j}<\infty$, then
\begin{equation}\label{1.1}
\lambda^{(a)}_{n}(\ell_{2},\gamma,\ell_{q})_{q}=\bigg(\int_{\ell_{2}}\|x-S_{n}x\|^{q}_{q}d\gamma(x)\bigg)^{\frac{1}{q}}=C(q)\Bigg(\sum\limits^{\infty}_{j=n+1}\sigma^{\frac{q}{2}}_{j}\Bigg)^{\frac{1}{q}}.
\end{equation}

{\rm(ii)} If $2<q<\infty$, then
\begin{equation}\label{1.2}
\lambda^{(a)}_{n}(\ell_{2},\gamma,\ell_{q})_{q}=C(q)\Bigg[(\kappa_{0}-n)^{\frac{q}{2}}\Bigg(\sum\limits^{\kappa_{0}}_{j=1}
\sigma^{\frac{r}{2}}_{j}\Bigg)^{1-\frac{q}{2}}+\sum\limits^{\infty}_{j=\kappa_{0}+1}\sigma^{\frac{q}{2}}_{j}\Bigg]^{\frac{1}{q}},
\end{equation}
where $\frac{1}{r}=\frac{1}{q}-\frac{1}{2}<0$, and $\kappa_{0}\equiv\kappa(n,q,C_{\gamma})<\infty$ is the biggest integer in $[n+1,\infty)$ such that
\begin{equation}\label{1.3}
(\kappa-n)\sigma^{\frac{r}{2}}_{\kappa}\leq\sum\limits^{\kappa}_{j=1}\sigma^{\frac{r}{2}}_{j}.
\end{equation}
\end{thm}

\begin{rem}
The authors in \cite{CNL, DW} studied  average linear widths
$\lambda^{(a)}_{n}(D:\ell_{2}^N\rightarrow\ell_{q}^N,\gamma)_{p}$
and probabilistic linear widths
$\lambda_{n,\delta}(D:\ell_{2}^N\rightarrow\ell_{q}^{N},\gamma)$
with a Gaussian measure $\gamma$ on $\ell_2^N$ for diagonal
matrices $D={\rm diag}(d_1,d_2,\dots, d_N)$ satisfying $d_1\ge
d_2\ge \dots\ge d_N \ge0$. Let $v=\gamma\circ D^{-1}$ be the
induced measure on $\ell_2^N$. Then  $v$ is also a Gaussian
measure on $\ell_2^N$ whose mean is zero and  correlation operator
$C_{v}=DC_\gamma D^{*}=DC_{\gamma}D$ has eigenvectors $e_{j}$ and
eigenvalues $v_{j}=d^{2}_{j}\sigma_{j}$, where $A^{*}$ is the conjugate operator of $A$. We have
$$
\lambda^{(a)}_{n}(D:\ell_{2}^N\rightarrow\ell_{q}^N,\gamma)_{p}=\lambda^{(a)}_{n}(\ell_{2}^N,v,
\ell_{q}^N)_{p} \ \ {\rm and}\ \
\lambda_{n,\delta}(D:\ell_{2}^N\rightarrow\ell_{q}^N,\gamma)=\lambda_{n,\delta}(\ell_{2}^N,v,
\ell_{q}^N).
$$

When $N=m<\infty$, $\sigma_j=1$, \cite[Theorems 2.1 and
4.1]{DW} give exact values of
$\lambda^{(a)}_{n}(D:\mathbb{R}^m\rightarrow\ell_{q}^m,\gamma_{m})_{q}$, $1\le q<\infty$ and asymptotic estimates of
$\lambda_{n,\delta}(D:\mathbb{R}^m\rightarrow\ell_{q}^m,\gamma_m),\,
1\le q\le2$ with the standard  Gaussian measure $\gamma_m$ on
$\Bbb R^m$ for diagonal matrices $D={\rm diag}(d_1,d_2,\dots,
d_m)$. Theorems \ref{thm1.1} and \ref{thm1.5} (i) reduce to
\cite[Theorems 2.1 and 4.1]{DW} when $v_j=d^{2}_{j}, \ 1\le j\le
m$ and $v_j=0,\ j\ge m+1$.
\end{rem}

According to \cite[Corollary 1]{WZ} or \cite[Theorem 1.2]{F}, we have
$$\lambda^{(a)}_{n}(\ell_{2},\gamma,\ell_{q})_{p}\asymp\lambda^{(a)}_{n}(\ell_{2},\gamma,\ell_{q})_{q}$$
for all $1\leq q\le \infty$ and $0<p<\infty$, where the constants
of equivalence depend only on $p$ and $q$. Here,  $A\asymp B$
means $A\ll B$ and $A\gg B$, and $A\ll B\ (A\gg B)$ means that
there exists an inessential positive constant $c$ such that $A\le
cB\ (A\ge cB)$. Thus, asymptotic estimates of
$\lambda^{(a)}_{n}(\ell_{2},\gamma,\ell_{q})_{p}$ for all $1\leq
q<\infty$ and $0<p<\infty$ follow immediately from Theorem
\ref{thm1.1}.

\begin{cor}\label{cor1.3}
Let $J_{n}(C_{\gamma})_{q}$ denote the expression on the right-hand side of (\ref{1.1}) and (\ref{1.2}) for $1\leq q\leq2$ and $2<q<\infty$, respectively. Then, for all $0<p<\infty$ and $1\leq q<\infty$,
\begin{equation}\label{1.4}
\lambda^{(a)}_{n}(\ell_{2},\gamma,\ell_{q})_{p}\asymp
J_{n}(C_{\gamma})_{q}.
\end{equation}
\end{cor}

Corollary  \ref{cor1.3} provides the asymptotic estimates of
$\lambda^{(a)}_{n}(\ell_{2},\gamma,\ell_{q})_{p}$ for all $1\leq
q<\infty$ and $0<p<\infty$ for  general $\{\sigma_j\}$. However,
it seems to be impossible to obtain such estimates for the case $q=\infty$. For
$\{\sigma_j\}$ satisfying certain conditions,  our second theorem
provides the asymptotic estimates of
$\lambda^{(a)}_{n}(\ell_{2},\gamma,\ell_{q})_{p}$ for $2< q\le
\infty$ and  $0<p<\infty$.

\begin{thm}\label{thm1.4}
Let $2<q\le\infty$ and  $0<p<\infty$.

{\rm(i)} If $\sigma_{2j}\asymp\sigma_{j}$ and $\{j\sigma_{j}\}$ is
nonincreasing, then
\begin{equation}\label{1.5}
\lambda^{(a)}_{n}(\ell_{2},\gamma,\ell_{q})_{p}\asymp
\lambda^{(a)}_{n}(\ell_{2},\gamma,\ell_{q})_{1}\asymp
\int_{\ell_2}\|x-S_nx\|_q d\gamma(x)\asymp
\left\{
\begin{aligned}
 &\sigma_{n+1}^{1/2}n^{1/q},&&2<q<\infty,\\
 &\sigma_{n+1}^{1/2}\sqrt{\ln n},&&q=\infty.
\end{aligned}
\right.
\end{equation}

{\rm(ii)} If
$\limsup\limits_{j\rightarrow\infty}\frac{\sigma_{j+1}}{\sigma_j}=\omega<1$,
then
\begin{equation}\label{1.6}
\lambda^{(a)}_{n}(\ell_{2},\gamma,\ell_{q})_{p}\asymp
\lambda^{(a)}_{n}(\ell_{2},\gamma,\ell_{q})_{1}\asymp
\int_{\ell_2}\|x-S_nx\|_q d\gamma(x)\asymp \sigma^{1/2}_{n+1}.
\end{equation}
\end{thm}

In  both cases,  we obtain from Theorem \ref{thm1.4} that the
Fourier partial summation operators $S_n$ are asymptotically
optimal in the average case setting for $2<q\le \infty$ and
$\{\sigma_j\}$ satisfying the above conditions.  We recall that
$\sum_{j=1}^\infty \sz_j<+\infty$. This means that
$\lim\limits_{j\to\infty} j\sigma_j=0$. Hence, the condition that
$\{j\sigma_{j}\}$ is nonincreasing is a mild condition.
\vspace{2mm}

Our third theorem provides  asymptotic estimates  of the
probabilistic linear widths
$\lambda_{n,\delta}(\ell_{2},\gamma,\ell_{q})$ for $1\leq q\leq2$
with general $\{\sigma_j\}$ and $2<q\le\infty$ with $\{\sigma_j\}$
satisfying certain conditions.

\begin{thm}\label{thm1.5}
There is a constant  $\delta_{0}\in(0,1/2)$ depending only on $q$ such that
for $\delta\in(0,\delta_{0}]$, we have the following asymptotics.

{\rm (i)} Let  $1\leq q\leq2$ and
$\sum^{\infty}_{j=1}\sigma^{\frac{q}{2}}_{j}<\infty$. Then
\begin{equation*}
\lambda_{n,\delta}(\ell_{2},\gamma,\ell_{q})\asymp
\Bigg(\sum\limits^{\infty}_{j=n+1}\sigma^{\frac{q}{2}}_{j}\Bigg)^{1/q}+\sqrt{\ln(1/\delta)}\Bigg(\sum^{\infty}_{j=n+1}\sigma^{\frac{r}{2}}_{j}\Bigg)^{1/r},
\end{equation*}
where $\frac{1}{r}=\frac{1}{q}-\frac{1}{2}$,   the constants of
equivalence depend only on $q$.

{\rm (ii)} Let $2<q\le \infty$.

{\rm \qquad (a)} If $\sigma_{2n}\asymp\sigma_{n}$ and
$\{j\sigma_{j}\}$ is nonincreasing, then
\begin{equation*}
\lambda_{n,\delta}(\ell_{2},\gamma,\ell_{q})\asymp\left\{
\begin{aligned}
 &\sigma^{1/2}_{n+1}n^{1/q}+\sigma^{1/2}_{n+1}\sqrt{\ln(1/\delta)},&&2<q<\infty,\\
 &\sigma^{1/2}_{n+1}\sqrt{\ln n}+\sigma^{1/2}_{n+1}\sqrt{\ln(1/\delta)},&&q=\infty.
\end{aligned}
\right.
\end{equation*}

{\rm \qquad (b)} If
$\limsup\limits_{j\rightarrow\infty}\frac{\sigma_{j+1}}{\sigma_{j}}=\omega<1$,
then
\begin{equation*}
\lambda_{n,\delta}(\ell_{2},\gamma,\ell_{q})\asymp
\sigma^{1/2}_{n+1}+\sigma^{1/2}_{n+1}\sqrt{\ln(1/\delta)}.
\end{equation*}
\end{thm}

The following corollary follows immediately from Corollary
\ref{cor1.3} and Theorems \ref{thm1.4} and \ref{thm1.5}.

\begin{cor}
There is a constant $\delta_{0}\in(0,1/2)$ depending only on $q$ such that
for $\delta\in(0,\delta_{0}]$, we have the following asymptotics.

{\rm(i)} If $1\leq q\leq\infty$, $\sigma_{j}=j^{-\rho}$  and $\rho>\max(1,2/q)$, then
\begin{equation*}
\lambda^{(a)}_{n}(\ell_{2},\gamma,\ell_{q})_{p}\asymp\bigg\{\begin{aligned}&n^{-\frac{\rho}{2}+\frac{1}{q}},
&&1\le q<\infty,\\ & n^{-\frac{\rho}{2}}\sqrt{\ln n}, &&q=\infty,
\end{aligned}
\end{equation*} and
\begin{equation*}
\lambda_{n,\delta}(\ell_{2},\gamma,\ell_{q})\asymp
\left\{
\begin{aligned}
 &n^{-\frac{\rho}{2}+\frac{1}{q}}+n^{-\frac{\rho}{2}+(\frac{1}{q}-\frac{1}{2})_+}\sqrt{\ln(1/\delta)},&&1\leq q<\infty,\\
  &n^{-\frac{\rho}{2}}\sqrt{\ln n}+n^{-\frac{\rho}{2}}\sqrt{\ln(1/\delta)},&&q=\infty.
\end{aligned}
\right.
\end{equation*}

{\rm(ii)} If $1\leq q\leq\infty$, $\sigma_{j}=j^{-\rho}A^{-2j}$,
$\rho\ge 0$ and $A>1$, then
\begin{equation*}
\lambda^{(a)}_{n}(\ell_{2},\gamma,\ell_{q})_{p}\asymp
n^{-\frac{\rho}{2}}A^{-n}, \end{equation*}and
\begin{equation*}
\lambda_{n,\delta}(\ell_{2},\gamma,\ell_{q})\asymp n^{-\frac{\rho}{2}}A^{-n}(1+\sqrt{\ln(1/\delta)}).
\end{equation*}
Here, $a_+=\max\{a,0\}$.
\end{cor}

The paper is organized as follows. Preliminary results about the
classical linear widths
 and about
the Gaussian measure in $\ell_{2}$ are provided in Section
\ref{sec2}. In Section \ref{sec3}, we determine the exact values
of the average linear widths
$\lambda^{(a)}_{n}(\ell_{2},\gamma,\ell_{q})_{q}$ for $1\leq
q<\infty$. Section \ref{sec4} is devoted to the asymptotic
estimates of average linear $n$-width
$\lambda^{(a)}_{n}(\ell_{2},\gamma,\ell_{q})_{p}$ for $2<q\le
\infty$ and $0<p<\infty$ with $\{\sigma_j\}$ satisfying certain
conditions. Finally, in Section \ref{sec5}  we give the  proof of
Theorem \ref{thm1.5}.

\section{Preliminaries}\label{sec2}

 For the
classical linear widths, the following exact values of
$\lambda_{n}(D(B_2), \ell_{q})$ are known, where $D={\rm
diag}(d_1,d_2,\dots,d_n,\dots)$ is a diagonal operator with
$$d_1\ge d_2\ge\dots\ge d_n\ge \dots,$$ and
$$D(B_p)=\{Dx=(d_1x_1,d_2x_2,\dots,d_nx_n,\dots): x\in B_p\},\ 1\le p\le \infty.$$

(i) If $1\leq q\leq2$ and $\sum^{\infty}_{j=1}d_{j}^{r}<\infty$,
we have (see \cite[Theorem 11.11.4]{PI})
\begin{equation}\label{2.1}
\lambda_{n}(D(B_2),
\ell_{q})=\sup_{x\in
B_{2}}\|Dx-DS_{n}x\|_{q}=\bigg(\sum^{\infty}_{j=n+1}d^{r}_{j}\bigg)^{1/r},
\end{equation}
where $\frac{1}{r}=\frac{1}{q}-\frac{1}{2}$. In particular,
$$\lambda_{n}(D(B_2), \ell_{2})=\sup_{x\in
B_{2}}\|Dx-DS_{n}x\|_{2}=d_{n+1}.$$

(ii) If $q=\infty$, we have (see \cite[Theorem 11.11.7]{PI})
\begin{equation}\label{2.2}
\lambda_{n}(D(B_2),
\ell_{\infty})=\lambda_{n}(D(B_1), \ell_{2})=\inf_{\substack{V\subset \ell_2\\ {\rm dim} V\le
n}}\max\limits_{j\in\mathbb{N}}d_{j}\min\limits_{y\in
V}\|e_{j}-y\|_{2}=\max_{k>n}\bigg(\frac{k-n}{\sum^{k}_{j=1}d_j^{-2}}\bigg)^{1/2}.
\end{equation}

For the case $2<q<\infty$, the exact values of
$\lambda_{n}(D(B_{2}),\ell_{q})$ remain unknown. However, it is
easy to see that for $2< q\le \infty$,
\begin{equation}\label{2.3}
\lambda_{n}(D(B_{2}),\ell_{\infty})\leq\lambda_{n}(D(B_{2}),\ell_{q})
\leq\lambda_{n}(D(B_{2}),\ell_{2})=d_{n+1},
\end{equation}and
\begin{equation}\label{2.4}
\sup_{x\in B_2}\|Dx-DS_nx\|_q=d_{n+1}.
\end{equation}

If $d_{n}\asymp d_{2n}$, then
\begin{equation}\label{2.5}
\lambda_{n}(D(B_{2}),\ell_{\infty})=
\max_{k>n}\bigg(\frac{k-n}{\sum^{k}_{j=1}d_j^{-2}}\bigg)^{1/2}\geq\bigg(\frac{n}{\sum^{2n}_{j=1}d_j^{-2}}\bigg)^{1/2}\geq\bigg(\frac{n}{2nd^{-2}_{2n}}\bigg)^{1/2}\gg
d_{n+1}.
\end{equation}And if  $\limsup\limits_{j\rightarrow\infty}\frac{d_{j+1}}{d_{j}}=\omega<1$,
then there exists an $N_0$ such that $\frac{d_{n+1}}{d_n}\le \tau$
for $n\ge N_0$, where $\tau=\frac{\omega+1}{2}<1$.  Hence,
 for $m>n\ge1$,
\begin{equation}\label{2.6}
\frac{d_m}{d_n}\le \frac{d_1}{d_{N_0}} \tau^{m-n-N_0}\ll
\tau^{m-n}.
\end{equation}
It follows from \eqref{2.2} and \eqref{2.6} that
\begin{align}\label{2.7}
&\lambda_{n}(D(B_{2}),\ell_{\infty})=
\max_{k>n}\bigg(\frac{k-n}{\sum^{k}_{j=1}d_j^{-2}}\bigg)^{1/2}\ge\bigg(\frac{1}{\sum^{n+1}_{j=1}d_j^{-2}}\bigg)^{1/2}\notag\\
&=\frac{d_{n+1}}{\Big(\sum^{n+1}_{j=1}\frac{d^{2}_{n+1}}{d^{2}_j}\Big)^{1/2}}\gg\frac{d_{n+1}}{\Big(\sum^{n+1}_{j=1}\tau^{2(n+1-j)}\Big)^{1/2}}\gg
d_{n+1}.
\end{align}
In both cases, by \eqref{2.3},  \eqref{2.4},  \eqref{2.5} and
\eqref{2.7} we obtain that  for $2\le q\le \infty$,
\begin{equation}\label{2.8}
\lambda_{n}(D(B_2),\ell_{q})\asymp  d_{n+1}=\sup_{x\in
B_2}\|Dx-DS_nx\|_q .
\end{equation}

Let $\mu$ be a Gaussian measure on $\ell_2$ with mean zero and
correlation operator $C_{\mu}$. Then $C_\mu$ is a symmetric
nonnegative nuclear operator on $\ell_2$, and   for any
$x,y\in\ell_{2}$, we have (see \cite[pp. 48-49]{Bo})
\begin{equation}\label{2.9}
\langle C_{\mu}x,y\rangle=\int_{\ell_{2}}\langle x,v\rangle\langle
y,v\rangle d\mu(v).
\end{equation}

Denote by $(\ell_{2})^{*}$ the space of all continuous linear
functionals on $\ell_{2}$, and by $L_{2}(\mu)$ the usual $L_2$
space of $\mu$-measurable functionals $f$ on $\ell_{2}$ with finite
norm
$$\|f\|_{L_{2}(\mu)}:=\Big(\int_{\ell_{2}}|f(x)|^{2}d\mu(x)\Big)^{1/2}.$$
Then $(\ell_2)^{*}=\ell_2$ can be embedded into $L_{2}(\mu)$. Set
(see \cite[p. 44]{Bo})
\begin{equation}\label{2.10}
H(\mu):=\Bigg\{x\in\ell_{2}:\|x\|_{H(\mu)}:=\sup_{\substack{y\in(\ell_{2})^{*}=\ell_{2}\\
R_{\mu}(y)(y)\leq1}}|\langle x,y\rangle|<\infty\Bigg\},\end{equation}
where
$$R_{\mu}(x)(y):=\int_{\ell_{2}}\langle x,v\rangle\langle
y,v\rangle d\mu(v)=\langle C_{\mu}x,y\rangle,\ \ \
x,y\in\ell_{2}.$$ The space $H(\mu)$ is called the
Cameron-Martin space (or the reproducing kernel Hilbert space) of
$\mu$. From \cite[Lemma 2.4.1]{Bo} we know   $x\in H(\mu)$ if and
only if there exists a $z\in \ell_2$ such that
$$x=C_\mu^{1/2}z,$$
and we have
$$\|x\|_{H(\mu)}=\|C_\mu^{1/2}z\|_{H(\mu)}=\|z\|_2.$$Furthermore,
for any $x,y\in\ell_2$, $$ \langle C_\mu^{1/2} x, C_\mu^{1/2}
y\rangle_{H(\mu)}=\langle x,y\rangle.$$

 Let $\gamma $ be the Gaussian measure on $\ell_2$ with
mean zero and correlation operator $C_{\gamma}$, where
$C_\gamma={\rm diag}(\sz_1,\dots,\sz_n,\dots)$ satisfying
$$\sz_1\ge\sz_2\ge\dots\ge \sz_n\ge \dots \ge0\ \ {\rm and}\ \  \sum_{j=1}^\infty \sz_j<\infty.$$ It is easy to see from
\eqref{2.10} that the Cameron-Martin space $H(\gamma)$ of the
Gaussian measure $\gamma$ is
$$H(\gamma)=\Bigg\{x\in \ell_2: \|x\|_{H(\gamma)}=\|C_\gamma^{-1/2}x\|_2=\Big(\sum^{\infty}_{j=1}\sigma^{-1}_{j}x^{2}_{j}\Big)^{1/2}<\infty\Bigg\},$$
with inner product $$ \ \langle x,y\rangle_{H(\gamma)}=\langle
C_\gamma^{-1/2}x,C_\gamma^{-1/2}y\rangle
=\sum^{\infty}_{j=1}\sigma^{-1}_{j}x_{j}y_{j}.$$

For any fixed $z_{1},\ldots,z_{n}\in\ell_{2}$, the random vector $\{\langle x,z_{j}\rangle\}^{n}_{j=1}$
on the measurable space $(\ell_{2},\gamma)$ has a centered Gaussian distribution with covariance matrix
 $$(\langle C_{\gamma}z_{i},z_{j}\rangle)^{n}_{i,j=1}=(\langle C_{\gamma}z_{i},C_{\gamma}z_{j}\rangle_{H(\gamma)})^{n}_{i,j=1}.$$
 In particular, on the cylindrical subsets, the measure $\gamma$ has the following form: let $g_{j},\, j=1,2,\ldots,n$ be any orthonormal system in
 $\ell_2$, with $g_j\in H(\gamma)$ for $j=1,2,\dots,n$,
 and let $K$ be an arbitrary Borel subset of $\mathbb{R}^{n}$. Then $C_\gamma^{1/2}g_{j}, j=1,2,\ldots,n$ is an orthonormal system in $H(\gamma)$,
  and the measure of the cylindrical subset
$$G=\{y\in\ell_{2}:\{\langle y, C^{-1/2}_{\gamma}g_{j}\rangle\}^{n}_{j=1}\in K\}$$
is equal to
$$\gamma(G)=(2\pi)^{-\frac{n}{2}}\int_{K}\exp\Big(-\frac{1}{2}\sum^{n}_{j=1}u^{2}_{j}\Big)du_{1}\cdots du_{n}.$$

For $2<q\le \infty$, we have $\ell_2\subset \ell_q$. Then the
above Gaussian measure $\gamma$
 can be extended to a measure $\widetilde{\gamma}$ on the space $\ell_{q}$  as follows:
  for any Borel subset $A$ of $\ell_{q}$,
  $\widetilde{\gamma}(A):=\gamma(A\cap\ell_{2}).$  For $1\le
q<2$, the space $\ell_q$ is a Borel measurable subspace of
 $\ell_2$. We give the condition for which $\gamma(\ell_q)=1$. In
 this case, the measure $\gamma$ may be viewed as a Gaussian
 measure on $\ell_q$ with the same Cameron-Martin space.

\begin{prop}\label{prop2.1} Let $\gamma$ be a Gaussian measure in $\ell_{2}$ with mean zero and correlation operator $C_{\gamma}$.
Then $\gamma(\ell_{q})=1$ if and only if $\sum^{\infty}_{j=1}\sigma^{q/2}_{j}<\infty$.
\end{prop}
\begin{proof}
 If $\sum^{\infty}_{j=1}\sigma^{q/2}_{j}<\infty$,
then from the proof of Lemma \ref{lem3.3}, we have
$$\int_{\ell_{2}}\|x\|_{q}^qd\gamma(x)\asymp
\sum^{\infty}_{j=1}\sigma^{q/2}_{j}<\infty,$$ which deduces that
$\|x\|_{q}<\infty$ almost surely on $\ell_2$. It follows that
$$\gamma(\{x\in\ell_{2}:\|x\|_{q}=\infty\})=\gamma(\ell_{2}\backslash\ell_{q})=0.$$
This implies
$$\gamma(\ell_{q})=\gamma(\ell_2)-\gamma(\ell_2\backslash\ell_{q})=1.$$

On the other hand, if $\gamma(\ell_{q})=1$, then the measure
$\gamma$ may be viewed as a Gaussian
 measure on $\ell_q$, and hence
$$\int_{\ell_{q}}\|x\|_{q}^qd\gamma(x)<\infty.$$ We note that
$$\int_{\ell_{q}}\|x\|_{q}^qd\gamma(x)=\int_{\ell_{2}}\|x\|_{q}^qd\gamma(x)\asymp \sum^{\infty}_{j=1}\sigma^{q/2}_{j}.$$
It follows that $\sum^{\infty}_{j=1}\sigma^{q/2}_{j}<\infty$. This
completes the proof of Proposition \ref{prop2.1}.
\end{proof}

For more information about the Gaussian measure on Banach spaces,
see \cite{Bo, Kuo}.

\section{Exact values of $\lambda^{(a)}_{n}(\ell_{2},\gamma,\ell_{q})_{q}$ for $1\le q<\infty$}\label{sec3}

In order to prove Theorem \ref{thm1.1},  we need several lemmas.
The first lemma ensures that the index $\kappa_{0}$ in Theorem
\ref{thm1.1} (ii) is finite.

\begin{lem}\label{lem3.1}(See \cite[Lemma 2.3]{NN})
Let $\{\delta_{j}\}$ be a positive increasing sequence, and $\lim\limits_{j\rightarrow\infty}\delta_{j}=+\infty$. Let $n\in\mathbb{N}$ and $n^{*}$ be the largest integer $m>n$ such that
\begin{equation}\label{3.1}
(m-n)\delta_{m}\leq\sum^{m}_{j=1}\delta_{j}.
\end{equation}
Then $n^{*}$ is finite, and for all $ n<m\le n^*$ the inequality
(\ref{3.1}) holds.
\end{lem}

The following lemma plays an essential role in supporting the proof of
Theorem \ref{thm1.1}.

\begin{lem}(See \cite[Lemma 11.11.6]{PI})\label{lem3.2}
Let $\{\xi_{j}\}$ be a scalar sequence, and let $n$ be a nonnegative integer. Then there exists an $n$-dimensional linear subspace $V$ of $\ell_{2}$ such that $\xi_{j}=\|P_{V}e_{j}\|^{2}_{2}$ for all $j\geq1$, if and only if $\sum^{\infty}_{j=1}\xi_{j}=n$ and $0\leq\xi_{j}\leq1$, where $P_{V}:\ell_{2}\rightarrow V$ denotes the orthogonal projection onto $V$.
\end{lem}

We need to compute the expectation of $\|X-TX\|^{q}_{q}$ with
$1\leq q<\infty$ and obtain the following exact value of
$\mathbb{E}\|X-TX\|^{q}_{q}$, where $X$ is a Gaussian random
vector in $\ell_{2}$ with mean zero and correlation operator
$C_{\gamma}$, and $T:\ell_{2}\rightarrow\ell_{q}$ is a bounded
linear operator with rank at most $n$.
\begin{lem}\label{lem3.3}
Let $1\leq q<\infty$, $\sum^{\infty}_{j=1}\sigma^{q/2}_{j}<\infty$
and $T:\ell_{2}\rightarrow\ell_{q}$ be a bounded linear operator with rank
at most $n$. Then
\begin{equation}\label{3.2}
\mathbb{E}\|X-TX\|^{q}_{q}=\int_{\ell_{2}}\|x-Tx\|^{q}_{q}d\gamma(x)=C(q)^{q}\sum^{\infty}_{j=1}\|\sigma^{\frac{1}{2}}_{j}e_{j}-C^{\frac{1}{2}}_{\gamma}T^{*}e_{j}\|_{2}^{q},
\end{equation}
where  $C(q)=\pi^{-\frac{1}{2q}}
2^{\frac{1}{2}}\,\Gamma\left(\frac{q+1}{2}\right)^{\frac{1}{q}},\
C_\gamma^{1/2}={\rm
diag}(\sigma_1^{1/2},\sigma_2^{1/2},\dots)$.
In particular, if $T=S_{n}$, then
\begin{equation}\label{3.3}
\int_{\ell_{2}}\|x-S_{n}x\|^{q}_{q}d\gamma(x)=C(q)^{q}\sum\limits^{\infty}_{j=n+1}\sigma^{\frac{q}{2}}_{j}.
\end{equation}
\end{lem}
\begin{proof}
For each fixed $j\in\mathbb{N}$, we set
$$L_{j}(x):=\langle x-Tx,e_{j}\rangle=\langle x,(I-T^{*})e_{j}\rangle,$$
where $I$ is the identity operator on $\ell_2$.
Obviously, $L_{j}$ is a bounded linear
functional on $\ell_{2}$. Since the distribution $\gamma$ is a symmetric
Gaussian measure on $\ell_{2}$, we know that $L_{j}(x)$, as a random
variable on the measurable space $(\ell_{2},\gamma)$, obeys the
normal distribution $N(0,R_{j}^{2})$, where
\begin{align}\label{3.4}
R_{j}^{2}&=\int_{\ell_{2}}|L_{j}(x)|^{2}d\gamma(x)=\int_{\ell_{2}}|\langle
x,(I-T^{*})e_{j}\rangle|^{2}d\gamma(x)=\langle
(I-T^{*})e_{j}, C_\gamma (I-T^{*})e_{j} \rangle\notag\\
&=\|C_\gamma^{1/2}(I-T^{*})e_{j}\|_2^2
=\|\sigma^{\frac{1}{2}}_{j}e_{j}-C^{\frac{1}{2}}_{\gamma}T^{*}e_{j}\|_{2}^2
\end{align}
is the variance of $L_{j}(x)$.
Since $L_{j}(x)=R_{j}\,\eta_j$ with $\eta_j\sim N(0,1)$, we obtain
\begin{align}\label{3.5}
\int_{\ell_{2}}|\langle
x-Tx,e_{j}\rangle|^{q}d\gamma(x)=\int_{\ell_{2}}|L_{j}(x)|^{q}d\gamma(x)=R_{j}^{q}\,\mathbb{E}|\eta_j|^{q}=C(q)^{q}R_{j}^{q},
\end{align}
where
$$C(q)^{q}=\mathbb{E}|\eta_j|^{q}=(2\pi)^{-\frac{1}{2}}\int_{\mathbb{R}}|t|^{q}e^{-\frac{t^{2}}{2}}dt=\pi^{-\frac{1}{2}}
2^{\frac{q}{2}}\Gamma\left(\frac{q+1}{2}\right).$$ It then follows from the monotone
convergence theorem that
$$\int_{\ell_{2}}\|x-Tx\|^{q}_{q}d\gamma(x)=\int_{\ell_{2}}
\sum^{\infty}_{j=1}|\langle x-Tx,e_{j}\rangle|^{q}d\gamma(x)=\sum^{\infty}_{j=1}\int_{\ell_{2}}|\langle x-Tx,e_{j}\rangle|^{q}d\gamma(x),$$
which, together with \eqref{3.5} and \eqref{3.4}, implies that
\begin{align*}
\int_{\ell_{2}}\|x-Tx\|^{q}_{q}d\gamma(x)
=C(q)^{q}\sum^{\infty}_{j=1}R_j^q=C(q)^{q}\sum^{\infty}_{j=1}\|\sigma^{\frac{1}{2}}_{j}e_{j}-C^{\frac{1}{2}}_{\gamma}T^{*}e_{j}\|_{2}^{q}.
\end{align*}

Finally, to complete the proof, it remains to note that (\ref{3.3}) is deduced directly from (\ref{3.2}) due to the fact that $S^{*}_{n}e_{j}=S_{n}e_{j}=e_{j}$ if $1\leq j\leq n$ and $S^{*}_{n}e_{j}=S_{n}e_{j}=0$ otherwise.
\end{proof}

The following lemma  shows that the average linear $n$-widths have
lower bounds  in terms of the extremal value of real function
defined on a compact and convex subset of $\ell_1$.
\begin{lem}\label{lem3.4}
If $1\leq q<\infty$, then
\begin{equation}\label{3.6}
\lambda^{(a)}_{n}(\ell_{2},\gamma,\ell_{q})_{q}\geq C(q)I_{q}(n),
\end{equation}
where $I_{q}(n)$ is defined by
\begin{equation}\label{3.7}
I_{q}(n):=\min\Bigg\{\bigg(\sum^{\infty}_{j=1}
\sigma^{\frac{q}{2}}_{j}(1-x_{j})^{\frac{q}{2}}\bigg)^{\frac{1}{q}}:
x=\{x_{j}\}\in[0,1]^{\infty}\quad{\rm
and}\quad\sum^{\infty}_{j=1}x_{j}=n\Bigg\}.
\end{equation}
\end{lem}

\begin{proof} Without loss of generality we assume that
$\sigma_j>0$ for all $j\in\Bbb N$. We start with the proof of the
lower estimate
$\lambda^{(a)}_{n}(\ell_{2},\gamma,\ell_{q})_{q}\geq
C(q)I_{q}(n)$. To this end, let
$T_{n}:\ell_{2}\rightarrow\ell_{q}$ be an arbitrary linear
operator with rank at most $n$, and let $\tilde
V=C_\gamma^{1/2}T^{*}_{n}(\ell_{2})$ be the range of the  operator
$C_\gamma^{1/2}T^{*}_{n}$. Then $\tilde V$ is the subspace of
$H(\gamma)$, and $\dim\tilde V\leq n$. Let $V$ be an
$n$-dimensional linear subspace of $H(\gamma)$ containing $\tilde
V$. We have
\begin{equation}\label{3.8}
\sum^{\infty}_{j=1}\|\sigma^{\frac{1}{2}}_{j}e_{j}-C^{\frac{1}{2}}_{\gamma}T_{n}^{*}e_{j}\|_{2}^{q}\geq\sum^{\infty}_{j=1}\sigma^{\frac{q}{2}}_{j}\min_{y\in V}\|e_{j}-y\|^{q}_{2}=\sum^{\infty}_{j=1}\sigma^{\frac{q}{2}}_{j}(1-\|P_{V}e_{j}\|^{2}_{2})^{\frac{q}{2}},
\end{equation}
where $P_{V}$ denotes the orthogonal projection onto $V$ in
$\ell_{2}$. Let
$$x_{j}=\|P_{V}e_{j}\|^{2}_{2},\quad{\rm for}\ j\geq1.$$
Obviously,  $0\leq x_{j}\leq1$ for each $j\geq1$, and
$$\sum^{\infty}_{j=1}x_{j}=\sum^{\infty}_{j=1}\|P_{V}e_{j}\|^{2}_{2}=\sum^{\infty}_{j=1}\langle
P_{V}e_{j},e_{j}\rangle={\rm trace}\ P_{V}= n.$$
Hence, by
\eqref{3.8} we have
$$
\sum^{\infty}_{j=1}\|\sigma^{\frac{1}{2}}_{j}e_{j}-C^{\frac{1}{2}}_{\gamma}T_{n}^{*}e_{j}\|_{2}^{q}
\ge\sum^{\infty}_{j=1}\sigma^{\frac{q}{2}}_{j}(1-x_{j})^{\frac{q}{2}}\geq I_{q}(n)^{q},
$$
which, together with \eqref{3.2}, yields that
$$\Big(\int_{\ell_{2}}\|x-T_{n}x\|^{q}_{q}d\gamma(x)\Big)^{\frac{1}{q}}\geq C(q)I_{q}(n).$$
Since $T_{n}$ is an arbitrary linear operator with rank at most $n$, the desired lower estimate is deduced immediately.
\end{proof}

\begin{lem}\label{lem3.5}
Let $1\leq q<\infty$. Then
\begin{equation}\label{3.9}
I_{q}(n)^{q}=\left\{
\begin{aligned}
 &\sum^{\infty}_{j=n+1}\sigma^{\frac{q}{2}}_{j},&&1\leq q\leq2,\\
 &(\kappa_{0}-n)^{\frac{q}{2}}\bigg(\sum^{\kappa_{0}}_{j=1}\sigma^{\frac{r}{2}}_{j}\bigg)^{1-\frac{q}{2}}+\sum^{\infty}_{j=\kappa_{0}+1}\sigma^{\frac{q}{2}}_{j},&&2<q<\infty,
\end{aligned}
\right.
\end{equation}
where $\frac{1}{r}=\frac{1}{q}-\frac{1}{2}$, and
$\kappa_{0}\equiv\kappa(n,q,C_{\gamma})$ is defined  in Theorem
\ref{thm1.1}.
\end{lem}

\begin{proof}
If $1\leq q\leq 2$, $\{x_{j}\}\in[0,1]^{\infty}$,
and $\sum^{\infty}_{j=1}x_{j}=n$, then from the monotonicity of
 $\{\sigma_{j}\}$, we have
\begin{align*}
\sum^{\infty}_{j=1}\sigma^{\frac{q}{2}}_{j}(1-x_{j})^{\frac{q}{2}}&\geq\sum^{\infty}_{j=1}\sigma^{\frac{q}{2}}_{j}(1-x_{j})
=\sum^{n}_{j=1}\sigma^{\frac{q}{2}}_{j}(1-x_{j})+\sum^{\infty}_{j=n+1}\sigma^{\frac{q}{2}}_{j}(1-x_{j})\\
&\geq\sigma^{\frac{q}{2}}_{n}\sum^{n}_{j=1}(1-x_{j})+\sum^{\infty}_{j=n+1}\sigma^{\frac{q}{2}}_{j}-\sigma^{\frac{q}{2}}_{n}\sum^{\infty}_{j=n+1}x_{j}\\
&=\sigma^{\frac{q}{2}}_{n}\bigg(n-\sum^{\infty}_{j=1}x_{j}\bigg)+\sum^{\infty}_{j=n+1}\sigma^{\frac{q}{2}}_{j}=\sum^{\infty}_{j=n+1}\sigma^{\frac{q}{2}}_{j},
\end{align*}
where the equalities can be achieved whenever
$x_{1}=\cdots=x_{n}=1$ and $x_{j}=0$ for $j\geq n+1$. Equality
\eqref{3.9} for $1\leq q\leq 2$ then follows.

Now it is sufficient to prove \eqref{3.9} for $2<q<\infty$. Let
$\{y_{j}\}\in[0,1]^{\infty}$ be the minimizer of
$$\sum^{\infty}_{j=1}\sigma^{q/2}_{j}(1-x_{j})^{q/2}$$ subject to
the conditions $\{x_{j}\}\in[0,1]^{\infty} \ {\rm
and}\  \sum^{\infty}_{j=1}x_{j}=n.$ Note that
$$\sum^{\infty}_{j=1}\sigma^{q/2}_{j}(1-x_{j})^{q/2}\geq\sum^{\infty}_{j=1}\sigma^{q/2}_{j}(1-x^{*}_{j})^{q/2},$$
where $\{x^{*}_{j}\}$ is a nonincreasing
rearrangement of  $\{x_{j}\}$. Therefore, without
loss of generality, we  assume that $1\ge y_{1}\geq
y_{2}\geq\ldots\geq0$. We will divide the proof into the following
four steps.\vspace{2mm}

(i) We show that $y_{1}<1$. If $y_{1}=1$, then there exists a
$j_{0}\in\mathbb{N}$ such that $y_{1}=\cdots=y_{j_{0}}=1,\
y_{j_{0}+1}<1.$ We set
$$\varphi(t)=\sigma^{q/2}_{j_{0}}(1-(y_{j_{0}}-t))^{q/2}+\sigma_{j_{0}+1}^{q/2}(1-(y_{j_{0}+1}+t))^{q/2}=
\sigma^{q/2}_{j_{0}}t^{q/2}+\sigma^{q/2}_{j_{0}+1}(1-y_{j_{0}+1}-t)^{q/2},$$
then for any $0\leq t\leq\delta$, $\varphi(t)$ achieves its
minimum at $t=0$, whenever $\delta>0$ is sufficiently small.
However,  since $q>2$ and
$$\varphi'(t)=\frac{q}{2}\sigma^{q/2}_{j_{0}}t^{q/2-1}-\frac{q}{2}\sigma^{q/2}_{j_{0}+1}(1-t-y_{j_{0}+1})^{q/2-1}<0$$
for $t>0$ sufficiently small, we have $\varphi(0)>\varphi(t)$
for $t\in(0,\delta]$. This leads to a contradiction. Hence
$y_{1}<1$.\vspace{2mm}

(ii) We show there exists a $j_0\in\mathbb{N}$ such that
$y_{j_0}=0$. Otherwise, for any $m\in\mathbb{N}$, $y_{m}>0$. We
set
$$\varphi_{m}(t)=
\sigma^{q/2}_{1}(1-(y_{1}+t))^{q/2}+\sigma_{m}^{q/2}(1-(y_{m}-t))^{q/2},
\ \ t\in[0,\delta_{0}],$$
where $\delta_{0}=\min\{y_{m},1-y_{1}\}.$ We
note that $\lim\limits_{m\rightarrow\infty}\sigma_{m}=0$. This
means that for $m$ sufficiently large,
$$\varphi'_{m}(t)=-\frac{q}{2}\sigma^{q/2}_{1}(1-(y_{1}+t))^{q/2-1}+\frac{q}{2}\sigma^{q/2}_{m}(1-y_{m}+t)^{q/2-1}<0$$
for $t\in(0,\delta]$ whenever $\delta$ is sufficiently small. This implies
$\varphi_{m}(t)< \varphi_{m}(0)$, leading to a
contradiction. Hence, there exists a $j_0\in \mathbb{N}$ such that
$y_{j_0}=0$.\vspace{2mm}

(iii) Choose $m$ such that $m>\max\{\kappa_0,j_0\}$. Then
$y_m=y_{m+1}=\dots=0$. Next, we show
\begin{equation}\label{3.10}
I_{q}(n)^{q}=\min\Bigg\{\sum^{m}_{j=1}\sigma^{\frac{q}{2}}_{j}x_{j}^{\frac{q}{2}}:
\{x_{j}\}^{m}_{j=1}\in[0,1]^{m}\quad{\rm and}\quad\sum^{m}_{j=1}x_{j}=m-n
\Bigg\}+\sum^{\infty}_{j=m+1}\sigma^{\frac{q}{2}}_{j}.
\end{equation}
It is easily seen that
\begin{align*}
I_{q}(n)^{q}&\leq\min\Bigg\{\sum^{m}_{j=1}\sigma^{\frac{q}{2}}_{j}(1-x_{j})^{\frac{q}{2}}+\sum^{\infty}_{j=m+1}\sigma^{\frac{q}{2}}_{j}:
\{x_{j}\}^{m}_{j=1}\in[0,1]^{m}\quad{\rm and}\quad\sum^{m}_{j=1}x_{j}=n
\Bigg\}\\
&=\min\Bigg\{\sum^{m}_{j=1}\sigma^{\frac{q}{2}}_{j}x_{j}^{\frac{q}{2}}:
\{x_{j}\}^{m}_{j=1}\in[0,1]^{m}\quad{\rm and}\quad\sum^{m}_{j=1}x_{j}=m-n\Bigg\}+\sum^{\infty}_{j=m+1}\sigma^{\frac{q}{2}}_{j}.
\end{align*}
On the other hand, since
$\{y_{j}\}\in[0,1]^{\infty}$ is the minimizer of $$\sum^{\infty}_{j=1}\sigma^{q/2}_{j}(1-x_{j})^{q/2}$$
  subject to the conditions $\{x_{j}\}\in[0,1]^{\infty} \ {\rm and}\  \sum^{\infty}_{j=1}x_{j}=n,$ it follows that

\begin{align*}I_{q}(n)^{q}&=\sum^{\infty}_{j=1}\sigma^{\frac{q}{2}}_{j}(1-y_{j})^{\frac{q}{2}}=\sum^{m}_{j=1}\sigma^{\frac{q}{2}}_{j}(1-y_{j})^{\frac{q}{2}}+\sum^{\infty}_{j=m+1}\sigma^{\frac{q}{2}}_{j}\\
&\geq\min\Bigg\{\sum^{m}_{j=1}\sigma^{\frac{q}{2}}_{j}(1-x_{j})^{\frac{q}{2}}:
\{x_{j}\}^{m}_{j=1}\in[0,1]^{m}\quad{\rm and}\quad\sum^{m}_{j=1}x_{j}=n\Bigg\}+\sum^{\infty}_{j=m+1}\sigma^{\frac{q}{2}}_{j}\\
&=\min\Bigg\{\sum^{m}_{j=1}\sigma^{\frac{q}{2}}_{j}x_{j}^{\frac{q}{2}}:
\{x_{j}\}^{m}_{j=1}\in[0,1]^{m}\quad{\rm and}\quad\sum^{m}_{j=1}x_{j}=m-n\Bigg\}+\sum^{\infty}_{j=m+1}\sigma^{\frac{q}{2}}_{j}.
\end{align*}
Hence, the equality \eqref{3.10} holds.\vspace{2mm}

(iv) We show that for any $2<q<\infty$,
$$I_{q}(n)^{q}=
(\kappa_{0}-n)^{\frac{q}{2}}\bigg(\sum^{\kappa_{0}}_{j=1}\sigma^{\frac{r}{2}}_{j}\bigg)^{1-\frac{q}{2}}+\sum^{\infty}_{j=\kappa_{0}+1}\sigma^{\frac{q}{2}}_{j}.$$
By \eqref{3.10} and \cite[(2.15)]{DW}, we get
\begin{align*}
I_{q}(n)^{q}&=\min\Bigg\{\sum^{m}_{j=1}\sigma^{\frac{q}{2}}_{j}x_{j}^{\frac{q}{2}}:
\{x_{j}\}^{m}_{j=1}\in[0,1]^{m}\quad{\rm and}\quad\sum^{m}_{j=1}x_{j}=m-n
\Bigg\}+\sum^{\infty}_{j=m+1}\sigma^{\frac{q}{2}}_{j}\\
&=(\kappa_{0,m}-n)^{\frac{q}{2}}\bigg(\sum^{\kappa_{0,m}}_{j=1}\sigma^{\frac{r}{2}}_{j}\bigg)^{1-\frac{q}{2}}+\sum^{m}_{j=\kappa_{0,m}+1}\sigma^{\frac{q}{2}}_{j}+\sum^{\infty}_{j=m+1}\sigma^{\frac{q}{2}}_{j}\\
&=(\kappa_{0,m}-n)^{\frac{q}{2}}\bigg(\sum^{\kappa_{0,m}}_{j=1}\sigma^{\frac{r}{2}}_{j}\bigg)^{1-\frac{q}{2}}+\sum^{\infty}_{j=\kappa_{0,m}+1}\sigma^{\frac{q}{2}}_{j},
\end{align*}
where
$\kappa_{0,m}\equiv\kappa(m,n,q,C_{\gamma})$ is the biggest
integer in $[n+1,m)$ satisfying \eqref{1.3}.
It follows from the definitions of $\kappa_{0,m}$ and $\kappa_{0}$
that $\kappa_{0,m}=\kappa_{0}$. Indeed, on the one hand, $\kappa_{0,m}$ is the biggest
integer in $[n+1,m)$ satisfying \eqref{1.3}
and since $m\geq\kappa_{0}$, we have $\kappa_{0,m}\geq\kappa_{0}$. On the other hand, since $\kappa_{0}$ is the biggest
integer in $[n+1,\infty)$ satisfying \eqref{1.3}, we have $\kappa_{0}\geq\kappa_{0,m}$. Thus, $\kappa_{0}=\kappa_{0,m}$. This implies that
$$I_{q}(n)^{q}=
(\kappa_{0}-n)^{\frac{q}{2}}\bigg(\sum^{\kappa_{0}}_{j=1}\sigma^{\frac{r}{2}}_{j}\bigg)^{1-\frac{q}{2}}+\sum^{\infty}_{j=\kappa_{0}+1}\sigma^{\frac{q}{2}}_{j}.$$
The desired estimate \eqref{3.9} for $2<q<\infty$ then follows.
\end{proof}

\begin{lem}\label{lem3.6}
For any $2<q<\infty$, we have
\begin{equation}\label{3.11}
\lambda^{(a)}_{n}(\ell_{2},\gamma,\ell_{q})_{q}\leq
C(q)\Bigg((\kappa_{0}-n)^{\frac{q}{2}}\bigg(\sum^{\kappa_{0}}_{j=1}\sigma^{\frac{r}{2}}_{j}\bigg)^{1-\frac{q}{2}}+
\sum^{\infty}_{j=\kappa_{0}+1}\sigma^{\frac{q}{2}}_{j}\Bigg)^{1/q}.
\end{equation}
\end{lem}

\begin{proof}
To prove the desired upper estimate, let
$$y_{k}=1-(\kappa_{0}-n)\sigma^{r/2}_{k}\bigg(\sum^{\kappa_{0}}_{j=1}\sigma^{r/2}_{j}\bigg)^{-1},\  \ k=1,2,\ldots,\kappa_{0},$$
and $y_{k}=0$ otherwise. Since $r=\frac{2q}{2-q}<0$ and $\sigma_{1}\geq\sigma_{2}\geq\cdots\geq \sigma_{n}\geq\cdots\geq0$, we deduce from the definition of $\kappa_{0}$ that $0\leq y_{k}\leq1$ for all $k\in\mathbb{N}$ and
$$\sum^{\infty}_{j=1}y_{j}=\sum_{j=1}^{\kappa_0}y_j=\kappa_0-(\kappa_0-n)=n.$$ Hence, by Lemma
\ref{lem3.2}, there exists a linear subspace $V\subset
S_{\kappa_{0}}(\ell_{2})\subset H(\gamma)$ with ${\rm dim}\,V= n$
and such that $y_{j}=\|P_{V}e_{j}\|^{2}_{2}$ for all
$j\in\mathbb{N}$, where $P_V$ is the orthogonal projection onto
$V$ in $\ell_2$.  Clearly, $\|e_{j}-P_{V}e_{j}\|^{2}_{2}=1-y_{j}$.
We remark that $P_VC_{\gamma}^{-1/2}$ is well defined on $\ell_2$
if $V\subset H(\gamma)$. Indeed, if $\{f_j\}_{j=1}^n$ is an
orthonormal basis in $V$ with respect to $\langle
\cdot,\cdot\rangle$, then for all $f_{j}$, $1\leq j\leq n$, there exists a unique $h_j\in\ell_{2}$, such that
$f_j=C_\gamma^{1/2}h_j$. For $x\in\ell_{2}$, we have
$$P_VC_\gamma^{-1/2}x=\sum_{j=1}^n\langle C_\gamma^{-1/2}x,
f_j\rangle f_j=\sum_{j=1}^n\langle C_\gamma^{-1/2}x,
C_\gamma^{1/2}h_j\rangle f_j=\sum_{j=1}^n\langle x, h_j\rangle
f_j.$$
Thus, using Lemma \ref{lem3.3}, we obtain
\begin{align*}
\lambda^{(a)}_{n}(\ell_{2},\gamma,\ell_{q})^{q}_{q}&\leq\int_{\ell_{2}}\|x-C^{1/2}_{\gamma}P_{V}C^{-1/2}_{\gamma}x\|^{q}_{q}d\gamma(x)
=C(q)^{q}\sum^{\infty}_{j=1}\sigma^{q/2}_{j}\|e_{j}-P_{V}e_{j}\|^{q}_{2}\\
&=C(q)^{q}\sum^{\infty}_{j=1}\sigma^{q/2}_{j}(1-y_{j})^{q/2}
=C(q)^{q}\Bigg((\kappa_{0}-n)^{\frac{q}{2}}\bigg(\sum^{\kappa_{0}}_{j=1}\sigma^{\frac{r}{2}}_{j}\bigg)^{1-\frac{q}{2}}
+\sum^{\infty}_{j=\kappa_{0}+1}\sigma^{\frac{q}{2}}_{j}\Bigg),
\end{align*}
which gives the desired upper estimates, and hence completes the proof of Lemma \ref{lem3.6}.
\end{proof}

Now we are in a position to prove Theorem \ref{thm1.1}.

\

\noindent{\it Proof of Theorem \ref{thm1.1}.}

By Lemmas \ref{lem3.3} and  \ref{lem3.6}, we obtain the upper
bound of $\lambda^{(a)}_{n}(\ell_{2},\gamma,\ell_{q})_{q}$ for
$1\le q\le2$ and $2<q<\infty$, respectively. By Lemmas
\ref{lem3.4} and \ref{lem3.5}, we get the lower bound
 of $\lambda^{(a)}_{n}(\ell_{2},\gamma,\ell_{q})_{q}$ for $1\le
 q<\infty$.
 The proof of Theorem \ref{thm1.1} is completed.\qed

\section{Asymptotic estimates of  $\lambda^{(a)}_{n}(\ell_{2},\gamma,\ell_{\infty})_{p}$ for $0<p<\infty$}\label{sec4}

We have obtained the exact values of
$\lambda^{(a)}_{n}(\ell_{2},\gamma,\ell_{q})_{q}$ for $1\le
q<\infty$ in Section \ref{sec3}. However, the exact values of
$\lambda^{(a)}_{n}(\ell_{2},\gamma,\ell_{\infty})_{p}$ remain
unknown for any $0<p<\infty$. From our proof of Theorem
\ref{thm1.1}, we know that there exists an $n$-dimensional linear
subspace $V$ of $H(\gamma)$ such that
\begin{equation}\label{4.1}
\lambda^{(a)}_{n}(\ell_{2},\gamma,\ell_{q})_{q}=\left(\int_{\ell_{2}}\|x-C_\gamma^{1/2}P_{V}C_\gamma^{-1/2}x\|^{q}_{q}d\gamma(x)\right)^{1/q},\quad
1\leq q<\infty,
\end{equation}
where $P_V$ is the orthogonal projection  onto
$V$ with respect to the inner product $\langle\cdot,\cdot \rangle
$. We remark that $C_\gamma^{1/2}P_{V}C_\gamma^{-1/2}$ is just the
orthogonal projection  onto
 the range $C_\gamma^{1/2}(V)$ of
$C_\gamma^{1/2}P_{V}C_\gamma^{-1/2}$ with respect to the inner
product $\langle\cdot,\cdot \rangle_{H(\gamma)}$. Indeed,
 if $\{f_j\}_{j=1}^n$ is an
orthonormal basis in $V$ with respect to $\langle
\cdot,\cdot\rangle$, then  $\{C_{\gamma}^{1/2}f_j\}_{j=1}^n$ is an
orthonormal basis in $C_{\gamma}^{1/2}(V)$ with respect to
$\langle \cdot,\cdot\rangle_{H(\gamma)}$, since
$$\langle C_{\gamma}^{1/2}f_i,
C_{\gamma}^{1/2}f_j \rangle_{H(\gamma)}=\langle
f_i,f_j\rangle=\delta_{ij},\ \ 1\le i,j\le n.$$ Thus, we obtain  for
$x\in H(\gamma)$,
$$C_{\gamma}^{1/2}P_VC_\gamma^{-1/2}x=\sum_{j=1}^n\langle C_\gamma^{-1/2}x,
f_j\rangle C_{\gamma}^{1/2} f_j=\sum_{j=1}^n\langle x,
C_{\gamma}^{1/2}f_j\rangle_{H(\gamma)} C_{\gamma}^{1/2}f_j,$$
which means that $C_{\gamma}^{1/2}P_VC_\gamma^{-1/2}$ is the orthogonal
projection onto
 $C_\gamma^{1/2}(V)$  with respect to  $\langle\cdot,\cdot \rangle_{H(\gamma)}$.

Our next theorem shows that \eqref{4.1} for $q=\infty$ is true for
some $n$-dimensional linear subspace $V$ of $H(\gamma)$.

\begin{thm}\label{thm4.1}
For $0<p<\infty$, we have
\begin{equation}\label{4.2}
\lambda^{(a)}_{n}(\ell_{2},\gamma,\ell_{\infty})^{p}_{p}=\inf\limits_{\substack{V\subset H(\gamma)\\
{\rm dim}V\le
n}}\int_{\ell_{2}}\|x-C_\gamma^{1/2}P_{V}C_\gamma^{-1/2}x\|^{p}_{\infty}d\gamma(x).
\end{equation}
\end{thm}

In order to prove Theorem \ref{thm4.1}, we need the following
 comparison lemma for
Gaussian measures (see \cite[p. 107]{Bo}).

\begin{lem}\label{lem4.2}
Let $X$ and $Y$ be the two Gaussian random vectors on $\ell_{2}$
with mean zero. If
\begin{equation}\label{4.3}
\mathbb{E}|\langle x,X\rangle|^{2}\leq\mathbb{E}|\langle
x,Y\rangle|^{2},\ \ \ \forall x\in\ell_{2},
\end{equation}
then for any $0<p<\infty$, and any semi-norm
$\varphi:\ell_{2}\rightarrow[0,\infty)$ on $\ell_{2}$,
$$\mathbb{E}\varphi(X)^{p}\leq\mathbb{E}\varphi(Y)^{p}.$$
\end{lem}
\begin{proof}
By   \cite[Theorem 3.3.6]{Bo},  (\ref{4.3}) is equivalent to the
following condition: For every absolutely convex Borel set
$A\in\ell_{2}$,
\begin{equation}\label{4.4}
\mathbb{P}(X\in A)\geq\mathbb{P}(Y\in A).
\end{equation}
We recall that a set $A$ is called absolutely convex if, for all
$x_{1},x_{2}\in A$ and $|\lambda_{1}|+|\lambda_{2}|\leq1$, we have
$\lambda_{1}x_{1}+\lambda_{2}x_{2}\in A$. Given any semi-norm
$\varphi$ in $\ell_{2}$, and any $t>0$, the set
$$A_{t}:=\{x\in\ell_{2}:\varphi(x)\leq t\}$$ is an absolutely convex
Borel set. Hence, by (\ref{4.4}) we obtain
$$\mathbb{P}(\varphi(X)>t)=1-\mathbb{P}(X\in A_{t})
\leq1-\mathbb{P}(Y\in A_{t})=\mathbb{P}(\varphi(Y)>t),$$ which
implies that
$$
\mathbb{E}\varphi(X)^{p}=p\int^{\infty}_{0}t^{p-1}\mathbb{P}(\varphi(X)>t)dt\leq
p\int^{\infty}_{0}t^{p-1}\mathbb{P}(\varphi(Y)>t)dt=\mathbb{E}\varphi(Y)^{p}.
$$
This completes the proof.
\end{proof}

\noindent {\it Proof of Theorem \ref{thm4.1}.}

Let $T:\ell_{2}\rightarrow\ell_{\infty}$ be a linear operator with
rank at most $n$, and let $V=C_\gamma ^{1/2}T^{*}(\ell_{2})$
denote the range of $C_\gamma ^{1/2}T^{*}$.
Obviously, $V$ is  a subspace of $H(\gamma)$, and $\dim V\le n$.
For the proof of Theorem \ref{thm4.1}, it suffices to prove that
for a Gaussian random vector $Z$ in $\ell_2$ with mean $0$ and
correlation operator $C_{\gamma}$ on $\ell_{2}$,
$$\mathbb{E}\|Z-TZ\|^{p}_{\infty}\geq\mathbb{E}\|Z-C_{\gamma}^{\frac{1}{2}}P_{V}C_{\gamma}^{-\frac{1}{2}}Z\|^{p}_{\infty}.$$
Since $Z-TZ$ and
$Z-C_{\gamma}^{\frac{1}{2}}P_{V}C_{\gamma}^{-\frac{1}{2}}Z$ are
both centered Gaussian random vectors on $\ell_{2}$, according to
Lemma \ref{lem4.2}, we reduce to proving that
$$
\mathbb{E}|\langle x,Z-TZ\rangle|^{2}\geq\mathbb{E}|\langle
x,Z-C_{\gamma}^{\frac{1}{2}}P_{V}C_{\gamma}^{-\frac{1}{2}}Z\rangle|^{2},\
\ \ x\in\ell_{2}.
$$
To achieve this, we utilize the formula
$$\mathbb{E}|\langle
x-T^{*}x,Z\rangle|^{2}=\|C_{\gamma}^{\frac{1}{2}}(I-T^{*})x\|^{2}_{2}.$$
Then, for all $x\in\ell_{2}$, we have
\begin{align*}
\mathbb{E}|\langle x,Z-TZ\rangle|^{2}&=\mathbb{E}|\langle x-T^{*}x,Z\rangle|^{2}=\|C_{\gamma}^{\frac{1}{2}}(I-T^{*})x\|^{2}_{2}
=\|C_{\gamma}^{\frac{1}{2}}x-C_{\gamma}^{\frac{1}{2}}T^{*}x\|^{2}_{2}\\
&\geq\|C_{\gamma}^{\frac{1}{2}}x-P_{V}C_{\gamma}^{\frac{1}{2}}x\|^{2}_{2}
=\mathbb{E}|\langle
x,Z-C_{\gamma}^{\frac{1}{2}}P_{V}C_{\gamma}^{-\frac{1}{2}}Z\rangle|^{2},
\end{align*}
which completes the proof  of Theorem \ref{thm4.1}.\qed

\vspace{2mm} Next, we are in a position to prove Theorem
\ref{thm1.4}. In order to prove \eqref{1.5} with $q=\infty$, we
need to use the result of Maiorov and Wasilkowski (see \cite[Theorem 2]{M3}).

\begin{lem}\label{lem4.3}
Let $m\ge 2n$ and let $\gamma_m$ be the standard Gaussian measure
on $\Bbb R^m$. Then
\begin{equation}\label{4.5}
\lambda^{(a)}_{n}(\mathbb{R}^{m},\gamma_{m},\ell_{\infty})_{1}\asymp\sqrt{\ln
(m-n)}.\end{equation}\end{lem} \vspace{2mm}

\noindent{\it Proof of Theorem \ref{thm1.4}.}

We start with the proof
of \eqref{1.5} for $2<q<\infty$. For the upper bound, we have
$$
\lambda^{(a)}_{n}(\ell_{2},\gamma,\ell_{q})_{q}
\leq\Big(\int_{\ell_{2}}\|x-S_{n}x\|^{q}_{q}d\gamma(x)\Big)^{1/q}
=C(q)\Bigg(\sum^{\infty}_{j=n+1}\sigma^{q/2}_{j}\Bigg)^{1/q}.
$$

We note that
\begin{align}\label{4.6}
\bigg(\sum^{\infty}_{j=n+1}\sigma^{q/2}_{j}\bigg)^{1/q}&=\bigg(\sum_{k=1}^\infty
\sum^{2^kn}_{j=2^{k-1}n+1}\sigma^{q/2}_{j}\bigg)^{1/q}
\le \bigg(\sum_{k=1}^{\infty} 2^{k-1}n\sigma^{q/2}_{2^{k-1}n}\bigg)^{1/q} \notag\\
&
=n^{1/q}\sigma_{n}^{1/2}\bigg(\sum^{\infty}_{k=1}2^{k-1}\Big(\frac{\sigma_{2^{k-1}n}}{\sigma_{n}}\Big)^{q/2}\bigg)^{1/q}\notag\\
&\leq cn^{1/q}\sigma_{n+1}^{1/2}\bigg(\sum^{\infty}_{k=1}
2^{-(k-1)(q/2-1)}\bigg)^{1/q}\leq Cn^{1/q}\sigma_{n+1}^{1/2},
\end{align}
where in the second inequality we used the fact that
$\{j\sigma_{j}\}$ is nonincreasing, and  the constant $C$ in the
last inequality is independent of $q$ for $q\geq q_{0}>2$. Hence, for $q_{0}\leq q<\infty$
\begin{equation}\label{4.7}
\lambda^{(a)}_{n}(\ell_{2},\gamma,\ell_{q})_{q}\leq CC(q)n^{1/q}\sigma_{n+1}^{1/2},
\end{equation}
and for $1\leq q<q_{0}$
$$\lambda^{(a)}_{n}(\ell_{2},\gamma,\ell_{q})_{q}\ll n^{1/q}\sigma_{n+1}^{1/2}.$$
For the lower bound, we note that if $\{j\sigma_{j}\}$ is nonincreasing, then
\begin{align*}
\sum^{\kappa-1}_{j=1}\sigma^{r/2}_{j}&=\sum^{\kappa-1}_{j=1}j^{r/2}\sigma^{r/2}_{j}j^{-r/2}
\leq\kappa^{r/2}\sigma^{r/2}_{\kappa}\sum^{\kappa-1}_{j=1}j^{-r/2}
\leq\kappa^{r/2}\sigma^{r/2}_{\kappa}\int^{\kappa}_{1}x^{-r/2}dx\\
&\leq\kappa^{r/2}\kappa^{-r/2+1}\sigma^{r/2}_{\kappa}\frac{2}{2-r}
=\kappa\sigma^{r/2}_{\kappa}\frac{2}{2-r}.
\end{align*}
Thus, if $(\kappa-n)\sigma^{r/2}_{\kappa}\leq\sum_{j=1}^{\kappa}\sigma^{r/2}_{j}$, then $\kappa\leq\frac{r-2}{r}(n+1)$ and $\frac{r-2}{r}=2-\frac{2}{q}\in[1,2)$, which implies $\kappa_{0}\leq \lfloor c\rfloor n$ and $c>1$, where $\lfloor c\rfloor$ denotes the biggest integer not greater than $c$.
From \eqref{1.2}, we know
\begin{equation}\label{4.8}
\lambda^{(a)}_{n}(\ell_{2},\gamma,\ell_{q})_{q}\geq C(q)\Bigg(\sum^{\infty}_{j=\lfloor c\rfloor n+1}\sigma^{q/2}_{j}\Bigg)^{1/q}\gg n^{1/q}\sigma^{1/2}_{n+1}.
\end{equation}
By \eqref{4.7}, \eqref{4.8} and \cite[Corollary 1]{WZ}, we have
$$\lambda^{(a)}_{n}(\ell_{2},\gamma,\ell_{q})_{p}\asymp\lambda^{(a)}_{n}(\ell_{2},\gamma,\ell_{q})_{q}\asymp n^{1/q}\sigma^{1/2}_{n+1}$$
for $2<q<\infty$ and $0<p<\infty$.

Next, we show \eqref{1.5} for $q=\infty$.
According to \cite[Corollary 1]{WZ}, it suffices to prove
that
$$\lambda^{(a)}_{n}(\ell_{2},\gamma,\ell_{\infty})_{1}\asymp
\sigma^{1/2}_{n+1}\sqrt{\ln n}.$$

For the upper bound, by Stirling's formula we obtain
$$\lim\limits_{q\rightarrow+\infty}\frac{C(q)}{q^{1/2}}=\lim\limits_{q\rightarrow+\infty}
\left(\frac{q-1}{q}\right)^{1/2}2^{\frac{1}{2q}}\exp\left(-\frac{q-1}{2q}\right)=e^{-1/2}.$$
Thus, taking $q=\ln n$, \eqref{4.7} implies
$$\lambda^{(a)}_{n}(\ell_{2},\gamma,\ell_{\infty})_{1}\leq\lambda^{(a)}_{n}(\ell_{2},\gamma,\ell_{q})_{1}
\leq\lambda^{(a)}_{n}(\ell_{2},\gamma,\ell_{q})_{q}\ll
\sigma^{1/2}_{n+1}\sqrt{\ln n},$$
which gives the desired upper estimate.

For the lower bound,  letting $T_n$ be an arbitrary linear operator with rank at most $n$,
we have
$$\int_{\ell_{2}}\|x-T_{n}x\|_{\infty}d\gamma(x)\geq
\int_{\ell_{2}}\|S_{2n}x-S_{2n}T_{n}x\|_{\infty}d\gamma(x).$$
We shall show
$$\int_{\ell_{2}}\|S_{2n}x-S_{2n}T_{n}x\|_{\infty}d\gamma(x)\geq
\int_{\ell_{2}}\|S_{2n}x-S_{2n}T_{n}S_{2n}x\|_{\infty}d\gamma(x).$$
Since $S_{2n}X-S_{2n}T_{n}X$ and $S_{2n}X-S_{2n}T_{n}S_{2n}X$ are
both centered Gaussian random vectors on $\ell_{2}$, using Lemma
\ref{lem4.2},  we reduce to showing that
$$\int_{\ell_{2}}|\langle S_{2n}x-S_{2n}T_{n}x,z\rangle|^{2} d\gamma(x)\geq
\int_{\ell_{2}}|\langle
S_{2n}x-S_{2n}T_{n}S_{2n}x,z\rangle|^{2} d\gamma(x),\quad
\forall z\in\ell_{2}.$$
To see this, we use the fact that
$$\int_{\ell_{2}}|\langle Tx,z\rangle|^{2}
d\gamma(x)=\int_{\ell_{2}}\langle x,T^*z\rangle  \langle
x,T^*z\rangle d\gamma(x)=\langle C_\gamma T^*z,T^*z\rangle=\langle
TC_\gamma T^*z,z\rangle$$ for all linear operator $T$.
We then get, for all $z\in\ell_{2}$,
\begin{align}\label{4.9}
\int_{\ell_{2}}|\langle S_{2n}x-S_{2n}T_{n}x,z\rangle|^{2}
d\gamma(x)=\langle
S_{2n}(I-T_{n})C_{\gamma}(I-T^{*}_{n})S_{2n}z,z\rangle
\end{align}
and
\begin{align}\label{4.10}
\int_{\ell_{2}}|\langle
S_{2n}x-S_{2n}T_{n}S_{2n}x,z\rangle|^{2} d\gamma(x)=\langle
S_{2n}(I-T_{n})S_{2n}C_{\gamma}S_{2n}(I-T^{*}_{n})S_{2n}z,z\rangle.
\end{align}
It is easy to see that the operator
$$S_{2n}(I-T_{n})(C_{\gamma}-S_{2n}C_{\gamma}S_{2n})(I-T^{*}_{n})S_{2n}$$
is nonnegative, i.e., for all $z\in \ell_2$, we have
\begin{align*}
&\quad\ \langle
S_{2n}(I-T_{n})(C_{\gamma}-S_{2n}C_{\gamma}S_{2n})(I-T^{*}_{n})S_{2n}z,z\rangle\\
&=\langle
(C_{\gamma}^{1/2}-S_{2n}C_{\gamma}^{1/2}S_{2n})(I-T^{*}_{n})S_{2n}z,
(C_{\gamma}^{1/2}-S_{2n}C_{\gamma}^{1/2}S_{2n})(I-T^{*}_{n})S_{2n}z\rangle
\ge 0.
\end{align*}
Hence, \eqref{4.9} and \eqref{4.10} imply
$$\int_{\ell_{2}}|\langle S_{2n}x-S_{2n}T_{n}x,z\rangle|^{2} d\gamma(x)\geq\int_{\ell_{2}}|\langle S_{2n}x-S_{2n}T_{n}S_{2n}x,z\rangle|^{2}d\gamma(x),$$
which deduces that
\begin{align*}
\int_{\ell_{2}}\|x-T_{n}x\|_{\infty}d\gamma(x)&\geq\int_{\ell_{2}}\|S_{2n}x-S_{2n}T_{n}x\|_{\infty}d\gamma(x)
\geq\int_{\ell_{2}}\|S_{2n}x-S_{2n}T_{n}S_{2n}x\|_{\infty}d\gamma(x)\\
&=\int_{\ell_{2}}\|S_{2n}x-S_{2n}T_{n}S_{2n}S_{2n}x\|_{\infty}d\gamma(x)=\int_{\Bbb
R^{2n}}\|y-S_{2n}T_{n}S_{2n}y\|_{\infty}d\mu(y)\\
&=:J,
\end{align*}
where $\mu$ is the Gaussian measure on $\Bbb R^{2n}$ with mean zero and
covariance matrix ${\rm diag}(\sz_{1},\dots,\sz_{2n})$. We denote by
$\gamma_{2n}$ the Gaussian measure on $\Bbb R^{2n}$ with mean zero
and covariance matrix $I_{2n}$, the $2n$ by $2n$ identity matrix.
Obviously, $S_{2n}T_{n}S_{2n}$ may be viewed as a linear
operator on $\mathbb R^{2n}=S_{2n}(\ell_2)$ with rank at most $n$,
and
$$\int_{\Bbb
R^{2n}}|\langle y-S_{2n}T_{n}S_{2n}y,z\rangle|^{2}d\mu(y)\ge \sz_{2n}\int_{\Bbb
R^{2n}}|\langle y-S_{2n}T_{n}S_{2n}y,z\rangle|^{2}d\gamma_{2n}(y)$$
for all $z\in\mathbb{R}^{2n}$.
Again, using Lemma \ref{lem4.2} and \eqref{4.5} we deduce
\begin{align}\label{4.11}
J\ge
\sz^{1/2}_{2n}\int_{\mathbb{R}^{2n}}\|y-S_{2n}T_{n}S_{2n}y\|_{\infty}d\gamma_{2n}(y)\ge
\sz^{1/2}_{2n}\lambda^{(a)}_{n}(\mathbb{R}^{2n},\gamma_{2n},\ell_{\infty})_{1}\gg
\sigma_{n+1}^{1/2}\sqrt{\ln n},
\end{align}
which gives the desired lower estimate.

Next, we are in a position to prove \eqref{1.6} for $2<q\leq\infty$. Since
\begin{equation}\label{4.12}
\lambda^{(a)}_{n}(\ell_{2},\gamma,\ell_{\infty})_{1}\leq\lambda^{(a)}_{n}(\ell_{2},\gamma,\ell_{q})_{1}
\ll\lambda^{(a)}_{n}(\ell_{2},\gamma,\ell_{2})_{2},
 \end{equation}
it suffices to prove the upper estimate of
$\lambda^{(a)}_{n}(\ell_{2},\gamma,\ell_{2})_{2}$
and the lower estimate of
$\lambda^{(a)}_{n}(\ell_{2},\gamma,\ell_{\infty})_{1}$.

Let us first establish the upper estimate. We note that $\limsup\limits_{j\to\infty}\frac{\sigma_{j+1}}{\sigma_{j}}=\omega<1$.
Then by \eqref{2.6},
we have the following estimate for the sum:
\begin{equation}\label{4.13}
\sum^{\infty}_{j=n+1}\sigma_{j}=\sigma_{n+1}\sum^{\infty}_{j=n+1}\frac{\sigma_{j}}{\sigma_{n+1}}\ll\sigma_{n+1}\sum^{\infty}_{j=n+1}\tau^{(j-n-1)}\ll\sigma_{n+1}.
\end{equation}
Therefore, using \eqref{1.1} we obtain
\begin{equation}\label{4.14}
\lambda^{(a)}_{n}(\ell_{2},\gamma,\ell_{2})_{2}=\Big(\sum\limits^{\infty}_{j=n+1}\sigma_{j}\Big)^{1/2}\ll\sigma^{1/2}_{n+1}.
\end{equation}

We now prove the lower estimate. Using \eqref{3.5} we have
\begin{align*}
\lambda^{(a)}_{n}(\ell_{2},\gamma,\ell_{\infty})_{1}
&=\inf_{T_{n}}\int_{\ell_{2}}\|x-T_{n}x\|_{\infty}d\gamma(x)
=\inf_{T_{n}}\int_{\ell_{2}}\max\limits_{j\in\mathbb{N}}|\langle (I-T^{*}_{n})e_{j},x\rangle| d\gamma(x)\\
&\geq\max\limits_{j\in\mathbb{N}}\inf_{T_{n}}\int_{\ell_{2}}|\langle (I-T^{*}_{n})e_{j},x\rangle |d\gamma(x)\asymp\max\limits_{j\in\mathbb{N}}\inf_{T_{n}}\|C^{1/2}_{\gamma}(I-T^{*}_{n})e_{j}\|_{2}.
\end{align*}
As in the proof of Lemma \ref{lem3.4}, let
$T_{n}:\ell_{2}\rightarrow\ell_{\infty}$ be an arbitrary bounded
linear operator with rank at most $n$, and let
$V=C^{1/2}_{\gamma}T^{*}_{n}(\ell_{2})$. Then $\dim V\leq n$.
We note that
$\limsup\limits_{j\to\infty}\frac{\sigma_{j+1}}{\sigma_{j}}=\omega<1$.
By \eqref{2.2} and \eqref{2.6}, we have
\begin{align}\label{4.15}
\lambda^{(a)}_{n}(\ell_{2},\gamma,\ell_{\infty})_{1}&\gg\max\limits_{j\in\mathbb{N}}\sigma^{1/2}_{j}\min\limits_{y\in V}\|e_{j}-y\|_{2}
\ge \lambda_{n}(C^{1/2}_{\gamma}(B_{1}),\ell_{2})=\max_{k>n}\left(\frac{k-n}{\sum^{k}_{j=1}\sigma^{-1}_{j}}\right)^{1/2}\notag\\
&\geq\left(\frac{1}{\sum^{n+1}_{j=1}\sigma^{-1}_{j}}\right)^{1/2}=\left(\frac{\sigma_{n+1}}{\sum^{n+1}_{j=1}\frac{\sigma_{n+1}}{\sigma_{j}}}\right)^{1/2}
\gg\left(\frac{\sigma_{n+1}}{\sum^{n+1}_{j=1}\tau^{(n+1-j)}}\right)^{1/2}\notag\\
&\gg\sigma^{1/2}_{n+1}.
\end{align}

A combination of  \eqref{4.12}, \eqref{4.14}, and \eqref{4.15} yields
$$\lambda^{(a)}_{n}(\ell_{2},\gamma,\ell_{q})_{p}\asymp\sigma^{1/2}_{n+1},$$
for $2<q\leq\infty$ and $0<p<\infty$,
which completes the proof of \eqref{1.6} for $2<q\leq\infty$. \qed

\vspace{2mm}

\section{Asymptotic estimates of the probabilistic linear widths}\label{sec5}

In order to prove Theorem \ref{thm1.5}, we need the following
 lemmas. The first is about the estimate of the
tail bound of a standard Gaussian distribution.
\begin{lem}\label{lem5.1}(See \cite[Lemma 1.1.3]{Bo}, \cite{M2})
Let $\eta\sim N(0,1)$. Then for all $t>0$
$$\sqrt{2}\pi^{-1/2}(t^{-1}-t^{-3})e^{-t^{2}/2}\leq\mathbb{P}(|\eta|\geq t)\leq\sqrt{2}\pi^{-1/2}t^{-1}e^{-t^{2}/2}.$$
\end{lem}

\noindent In particular, if $\delta\in(0,e^{-1}]$, then
$\mathbb{P}(|\eta|\geq c_{0}\sqrt{\ln\delta^{-1}})>\delta$
for some absolute constant $c_{0}>0$.

\begin{lem}\label{lem5.2}
(See \cite[(3.2)]{LT}) Let $1\leq q\leq\infty$, let $\frac{1}{q}+\frac{1}{q'}=1$, and let $Y$ be a centered
Gaussian random vector in $\ell_2$. Then for all $t>0$
\begin{equation}\label{5.1}
\mathbb{P}(|\|Y\|_{q}-\mathbb{E}\|Y\|_{q}|\geq t)\leq2e^{-\frac{t^{2}}{2\sigma^{2}}},
\end{equation}
where
$$\sigma=\sup\limits_{\xi\in B_{q'}}(\mathbb{E}|\langle \xi,Y\rangle|^{2})^{\frac{1}{2}}.$$
\end{lem}
Inequality (\ref{5.1}) is called the Gaussian concentration
inequality.

Now we are in a position to prove Theorem \ref{thm1.5}.

\vspace{2mm}

\noindent{\it Proof of Theorem \ref{thm1.5}.}

Firstly, we prove
\begin{equation}\label{5.2}
\lambda_{n,\delta}(\ell_{2},\gamma,\ell_{q})\ll e^{(a)}_{q}(S_{n})+e_{q}(S_{n})\sqrt{\ln(2/\delta)}
\end{equation}
for all $1\leq q\leq\infty$, where
$$e_{q}(S_{n}):=\sup\limits_{\xi\in B_{2}}\|C^{\frac{1}{2}}_{\gamma}(I-S_{n})\xi\|_{q}$$
and
$$e^{(a)}_{q}(S_{n}):=\int_{\ell_{2}}\|x-S_{n}x\|_{q}d\gamma(x).$$
Let $F(x)=\|x-S_{n}x\|_{q}$ for all $x\in\ell_{2}$. Applying
\cite[Lemma 2.2.2]{Bo}, we can derive that $Y=(I-S_{n})X$ is also
a centered Gaussian random vector in $\ell_{2}$, and the Gaussian
concentration inequality (\ref{5.1}) shows
$$\mathbb{P}(|F(X)-\mathbb{E}(F(X))|\geq t)\leq2e^{-\frac{t^{2}}{2\sigma^{2}}},\qquad \forall t>0,$$
where, here and in what follows, $X$ is a Gaussian random vector
in $\ell_{2}$ with mean $0$ and correlation operator $C_{\gamma}$. For $1\leq q\leq\infty$, we have, by duality,
\begin{align*}
\sigma&=\sup\limits_{\xi\in B_{q'}}(\mathbb{E}|\langle
\xi,Y\rangle|^{2})^{1/2}
=\sup\limits_{\xi\in B_{q'}}(\mathbb{E}|\langle(I-S_{n})\xi,X\rangle|^{2})^{1/2}
=\sup\limits_{\xi\in B_{q'}}\|C^{\frac{1}{2}}_{\gamma}(I-S_{n})\xi\|_{2}\\
&=\sup\limits_{\xi\in B_{2}}\|C^{\frac{1}{2}}_{\gamma}(I-S_{n})\xi\|_{q}=e_{q}(S_{n}).
\end{align*}
This means that for
$$Q_{\delta}=\{x\in\ell_{2}:F(x)>e^{(a)}_{q}(S_{n})+e_{q}(S_{n})\sqrt{2\ln(2/\delta)}\}$$
with $\delta\in(0,1)$,
$$\gamma(Q_{\delta})\leq\mathbb{P}(|F(X)-e^{(a)}_{q}(S_{n})|>e_{q}(S_{n})\sqrt{2\ln(2/\delta)})\leq\delta.$$
By the definition of the linear $(n,\delta)$-widths, this last inequality  deduces the desired upper estimate
\begin{equation}\label{5.3}
\lambda_{n,\delta}(\ell_{2},\gamma,\ell_{q})\leq\sup\limits_{\ell_{2}\backslash{Q_{\delta}}}\|x-S_{n}x\|_{q}\ll e^{(a)}_{q}(S_{n})+e_{q}(S_{n})\sqrt{\ln(2/\delta)}.
\end{equation}

Next, we show
\begin{equation}\label{5.4}
\lambda_{n,\delta}(\ell_{2},\gamma,\ell_{q})\gg \sigma_{2}+\sigma_{1}\sqrt{\ln(2/\delta)},
\end{equation}
where
$$\sigma_{1}:=\lambda_{n}(C_{\gamma}^{1/2}(B_{2}),\ell_{q})$$
and
$$\sigma_{2}:=\lambda^{(a)}_{n}(\ell_{2},\gamma,\ell_{q})_{1}.$$
Hence, we need only to prove both
\begin{equation}\label{5.5}
\lambda_{n,\delta}(\ell_{2},\gamma,\ell_{q})\geq c\sigma_{1}\sqrt{\ln(2/\delta)},
\end{equation}
and
\begin{equation}\label{5.6}
\lambda_{n,\delta}(\ell_{2},\gamma,\ell_{q})\geq c\sigma_{2},
\end{equation}
for some absolute constant $c>0$.

For the proof of (\ref{5.5}), it  suffices  to prove that for any
bounded linear operator $T_{n}:\ell_{2}\rightarrow\ell_{q}$ with
rank at most $n$,
\begin{equation}\label{5.7}
\mathbb{P}(\|X-T_{n}X\|_{q}\geq c\sigma_{1}\sqrt{\ln(2/\delta)})>\delta.
\end{equation}
Indeed, (\ref{5.7}) shows that the set $A=\{x\in\ell_{2}:\|x-T_{n}x\|_{q}\geq c\sigma_{1}\sqrt{\ln(2/\delta)}\}$ cannot be entirely contained in any Borel subset $G\subset\ell_{2}$ with $\gamma(G)\leq\delta$. Hence, $A\backslash G\ne\emptyset$, and
$$\lambda_{n,\delta}(\ell_{2},\gamma,\ell_{q})\geq\inf_{G}\inf_{T_{n}}\sup_{x\in A\backslash G}\|x-T_{n}x\|_{q}\geq c\sigma_{1}\sqrt{\ln(2/\delta)}.$$
Letting
$\langle X,(I-T^{*}_{n})y\rangle=\|C^{\frac{1}{2}}_{\gamma}(I-T^{*}_{n})y\|_{2}\eta$, $\eta\sim N(0,1)$, we observe that
\begin{align*}
\mathbb{P}(\|X-T_{n}X\|_{q}\geq c\sigma_{1}\sqrt{\ln(2/\delta)})&=\mathbb{P}(\sup\limits_{y\in B_{q'}}|\langle X,(I-T^{*}_{n})y\rangle|\geq c\sigma_{1}\sqrt{\ln(2/\delta)})\\
&\geq\sup\limits_{y\in B_{q'}}\mathbb{P}(|\langle X,(I-T^{*}_{n})y\rangle|\geq c\sigma_{1}\sqrt{\ln(2/\delta)})\\
&=\sup\limits_{y\in B_{q'}}\mathbb{P}\Bigg(|\eta|\geq c\frac{\sigma_{1}\sqrt{\ln(2/\delta)}}{\|C^{1/2}_{\gamma}(I-T^{*}_{n})y\|_{2}}\Bigg)\\
&\geq\mathbb{P}\Bigg(|\eta|\geq c\frac{\sigma_{1}\sqrt{\ln(2/\delta)}}{B}\Bigg),
\end{align*}
where
$$B:=\sup\limits_{y\in B_{q'}}\|C^{\frac{1}{2}}_{\gamma}(I-T^{*}_{n})y\|_{2}
=\sup\limits_{y\in B_{2}}\|(I-T_{n})C^{\frac{1}{2}}_{\gamma}y\|_{q}
\geq\lambda_{n}(C^{\frac{1}{2}}_{\gamma}(B_{2}),\ell_{q})=\sigma_{1}.$$
It then follows by Lemma \ref{lem5.1} that
$$\mathbb{P}(\|X-T_{n}X\|_{q}\geq c\sigma_{1}\sqrt{\ln(2/\delta)})\geq\mathbb{P}(|\eta|\geq c\sqrt{\ln(2/\delta)})>\delta$$
with $c=c_{0}$ being the same as in Lemma \ref{lem5.1}. This shows (\ref{5.7}), and hence the lower estimate (\ref{5.5}).

Now it remains to prove the lower estimate (\ref{5.6}). Again, it suffices to show that
\begin{equation}\label{5.8}
\mathbb{P}(\|X-T_{n}X\|_{q}\geq c\sigma_{2})>\delta
\end{equation} holds for any bounded  linear operator $T_{n}:\ell_{2}\rightarrow\ell_{q}$ with rank at most $n$. Let
\begin{equation}\label{5.9}
e(T_{n},\gamma,t):=\inf\{\rho>0:\mathbb{P}(\|X-T_{n}X\|_{q}>\rho)\leq
t\},\ \ \ 0\leq t\leq1.
\end{equation}
We reduce the proof of (\ref{5.8}) to showing
\begin{equation}\label{5.10}
e(T_{n},\gamma,\delta)>c\sigma_{2}.
\end{equation}

Since the function $t\rightarrow e(T_{n},\gamma,t)$ defined by (\ref{5.9}) is a nonincreasing rearrangement of the function $G(x):=\|x-T_{n}x\|_{q}$ with respect to the Gaussian measure $\gamma$ on $\ell_{2}$, we have
\begin{equation}\label{5.11}
\mathbb{E}G(X)=\int_{\ell_{2}}G(x)d\gamma(x)=\int^{1}_{0}e(T_{n},\gamma,t)dt.
\end{equation}
From \cite[Lemma 2.2.2]{Bo}, we conclude that
$(I-T_{n})X$ is a centered Gaussian random vector in $\ell_{2}$,
and the Gaussian concentration inequality (\ref{5.1}) implies
\begin{equation}\label{5.12}
\mathbb{P}(|G(X)-\mathbb{E}(G(X))|\geq s)\leq2e^{-\frac{s^{2}}{2\sigma^{2}}},\ \ \ s\geq0.
\end{equation}
For $1\leq q<\infty$, we use Lemma \ref{lem3.3} to give
\begin{align*}
\sigma&=\sup\limits_{\xi\in B_{q'}}(\mathbb{E}|\langle \xi,X-T_{n}X\rangle|^{2})^{\frac{1}{2}}=\sup\limits_{\xi\in B_{q'}}(\mathbb{E}|\langle(I-T^{*}_{n})\xi,X\rangle|^{2})^{\frac{1}{2}}=\sup\limits_{\xi\in B_{q'}}\|C^{\frac{1}{2}}_{\gamma}(I-T^{*}_{n})\xi\|_{2}\\
&=\sup\limits_{\xi\in B_{q'}}\Bigg\|\sum\limits^{\infty}_{j=1}\xi_{j}C^{\frac{1}{2}}_{\gamma}(I-T^{*}_{n})e_{j}\Bigg\|_{2}
\leq\sup\limits_{\xi\in B_{q'}}\sum\limits^{\infty}_{j=1}|\xi_{j}|\|C^{\frac{1}{2}}_{\gamma}(I-T^{*}_{n})e_{j}\|_{2}=\Bigg(\sum\limits^{\infty}_{j=1}\|C^{\frac{1}{2}}_{\gamma}(I-T^{*}_{n})e_{j}\|^{q}_{2}\Bigg)^{1/q}\\
&=C(q)^{-1}(\mathbb{E}\|(I-T_{n})X\|^{q}_{q})^{1/q}\leq
c\mathbb{E}G(X)=:\sigma_{3},
\end{align*}where $C(q)^q=\pi^{-1/2}2^{q/2}\Gamma((q+1)/2)$.
For $q=\infty$, we have
\begin{align*}
\sigma&=\sup\limits_{\xi\in B_{1}}(\mathbb{E}|\langle \xi,X-T_{n}X\rangle|^{2})^{1/2}=\sup\limits_{\xi\in B_{1}}\Bigg\|\sum\limits^{\infty}_{j=1}\xi_{j}C^{\frac{1}{2}}_{\gamma}(I-T^{*}_{n})e_{j}\Bigg\|_{2}\leq\sup\limits_{\xi\in B_{1}}\sum\limits^{\infty}_{j=1}|\xi_{j}|\|C^{\frac{1}{2}}_{\gamma}(I-T^{*}_{n})e_{j}\|_{2}\\
&\leq\max_{j\in\mathbb{N}}\|C^{\frac{1}{2}}_{\gamma}(I-T^{*}_{n})e_{j}\|_{2}
=\max_{j\in\mathbb{N}}C(1)^{-1}\mathbb{E}|\langle(I-T^{*}_{n})e_{j},X\rangle|
\leq C(1)^{-1}\mathbb{E}\|(I-T_{n})X\|_{\infty}\\
&=c\mathbb{E}G(X)=:\sigma_{3}.
\end{align*}

Given $t\in(0,1]$, we choose $s=\sigma_{3}\sqrt{2\ln(2/t)}=C\mathbb{E}G(X)\sqrt{\ln(2/t)}$ such that $2\exp(-s^{2}/2\sigma^{2}_{3})=t$, and then apply (\ref{5.12}) to obtain
$$\mathbb{P}(|G(X)|>(1+C\sqrt{\ln(2/t)})\mathbb{E}G(X))\leq t.$$
Using (\ref{5.9}), we conclude that
\begin{equation}\label{5.13}
e(T_{n},\gamma,t)\leq(1+C\sqrt{\ln(2/t)})\mathbb{E}G(X).
\end{equation}
Without loss of generality we assume $C\ge 1$. Now letting
$\delta_{0}\in(0,1)$ be such that
$\int^{\delta_{0}}_{0}(1+\sqrt{\ln(2/t)})dt=\frac{1}{2C}$, we
deduce from (\ref{5.11}) and (\ref{5.13}) that for any
$\delta\in(0,\delta_{0})$,
\begin{align*}
\mathbb{E}G(X)&=\int^{1}_{0}e(T_{n},\gamma,t)dt
\leq\bigg(\int^{\delta}_{0}(1+C\sqrt{\ln(2/t)})dt\bigg)\mathbb{E}G(X)+\int^{1}_{\delta}e(T_{n},\gamma,t)dt\\
&\leq\frac{1}{2}\mathbb{E}G(X)+e(T_{n},\gamma,\delta).
\end{align*}
This implies that
$$e(T_{n},\gamma,\delta)\geq\frac{1}{2}\mathbb{E}G(X)\geq c\inf_{T_{n}}\mathbb{E}\|X-T_{n}X\|_{q}= c\lambda^{(a)}_{n}(\ell_{2},\gamma,\ell_{q})_{1}=c\sigma_{2},
$$
which proves (\ref{5.10}), and hence the desired lower estimate (\ref{5.6}).

For $1\le q\le 2$, Theorem \ref{thm1.5} (i) follows from
\eqref{5.2}, \eqref{5.4}, \eqref{2.1}, and \eqref{1.1}
immediately.
 For $2<q\le \infty$, Theorem \ref{thm1.5} (ii) follows from \eqref{5.2}, \eqref{5.4}, \eqref{2.8},
  and Theorem \ref{thm1.4} immediately.

  Theorem \ref{thm1.5} is proved.
\qed

\

\noindent\textbf{Acknowledgments}
  The authors were
supported by the National Natural Science Foundation of China
(Project no. 12371098).


\begin{thebibliography}{99}
\bibitem{Be}J. Bergh, J. L\"{o}fstr\"{o}m, Interpolation Spaces: An Introduction. Grundlehren der mathematischen Wissenschaften 223, Springer-Verlag, 1976.

\bibitem{Bo}V.I. Bogachev, Gaussian Measures. Mathematical surveys and Monographs, 62. Providence: Amer. Math Soc, 1998.

\bibitem{CNL}G. Chen, P. Nie, X. Luo, The approximation characteristic of diagonal matrix in probabilistic setting, J. Complexity 26 (2010) 336-343.

\bibitem{DW}F. Dai, H. Wang, Linear $n$-widths of diagonal matrices in the average and probabilistic settings, J. Funct. Anal. 262 (2012) 4103-4119.

\bibitem{DD}D. D\~{u}ng, Weighted sampling recovery of functions with mixed smoothness, J. Complexity 93 (2026) 102000.

\bibitem{FY}G. Fang, P. Ye, Probabilistic and average linear widths of Sobolev space with Gaussian measure, J. Complexity 19 (2003) 73-84.

\bibitem{FY2004}G. Fang, P. Ye, Probabilistic and average widths of Sobolev space with Gaussian
measure in $L_{\infty}$-Norm, Constr. Approx. 20(1) (2004) 159-172.

\bibitem{F}X. Fernique, Integrabilite des vecteurs gaussiens, C. R. Acad. Sci., Paris A 270 (1970) 1698-1699.

\bibitem{Hi}F.J. Hickernell, H. Wozniakowski, Integration and approximation in arbitrary
dimensions, Adv. Comput. Math. 12 (2000) 25-58.

\bibitem{H}K. Hollig, Approximationszahlen von Sobolev-Einbettungen, Dissertation, Bonn, 1979.

\bibitem{Kuo}H.H. Kuo, Gaussian measure in Banach space, Lecture Notes in Mathematics, 463. Berlin: Springer, 1975.

\bibitem{La}F.M. Larkin, Optimal approximation in Hilbert spaces with reproducing kernel functions, Math. Comp. 24 (112) (1970) 911-921.

\bibitem{LT}M. Ledoux, M. Talagrand, Probability in Banach Spaces: Isoperimetry and processes. Springer, 2013.

\bibitem{LGM}G.G. Lorentz, M.V. Golitschek, Y. Makovoz, Constructive Approximation, Advanced Problems, Springer-Verlag, New York, 1996.

\bibitem{M1}V.E. Maiorov, Kolmogorov's $(n,\delta)$-widths of spaces of the smooth functions, Russian Acad. Sci. Sb. Math. 79 (2) (1994) 265-279.

\bibitem{M2}V.E. Maiorov, Linear widths of function spaces equipped with the Gaussian measure, J. Approx. Theory 77 (1) (1994) 74-88.

\bibitem{M3}V.E. Maiorov, G.W. Wasilkowski, Probabilistic and average linear widths in $L_{\infty}$-norm with respect to $r$-fold Wiener measure, J. Approx. Theory 84 (1) (1996) 31-40.

\bibitem{MR}Y.V. Malykhin, K.S. Ryutin, Widths and rigidity of unconditional sets and random vectors. (Russian. Russian summary)
Izv. Ross. Akad. Nauk Ser. Mat. 89(2) (2025) 45-59; translation in
Izv. Math. 89(2) (2025) 261-273.

\bibitem{NN}V.K. Nguyen, V.D. Nguyen, Best $n$-term approximation of diagonal operators and application to function spaces with mixed smoothness, Analysis Math. 48 (4) (2022) 1127-1152.

\bibitem{P1990}A. Papageorgiou, G.W. Wasilkowski, On the average complexity of multivariate
problems, J. Complexity 6(1) (1990) 1-23.

\bibitem{P1993}S.H. Paskov, Average case complexity of multivariate integration for smooth
 functions, J. Complexity 9(2) (1993) 219-312.

\bibitem{PI}A. Pietsch, Operator Ideals, North-Holland, Amsterdam, 1980.

\bibitem{P}A. Pinkus, $n$-widths in Approximation Theory, Springer, Berlin, 1985.

\bibitem{R}K. Ritter, Average-Case Analysis of Numerical Problems, in Lecture Notes in Math, vol. 1733, Springer-Verlag, 2000.

\bibitem{SS}M.S. Shabozov, M.S. Saidusainov, Widths of some classes of analytic functions in the weighted Bergman space, Lobachevskii J. Math. 46(7) (2025) 3577-3587.
\bibitem{S1959}A.V. Sul'din, Wiener measure and its applications to approximation theory I, Izv.
Vyssh. Uchebn. Zaved. Mat. No.13(1959), 145-158(Russian).

\bibitem{S1960}A.V. Sul'din, Wiener measure and its applications to approximation theory II, Izv.
Vyssh. Uchebn. Zaved. Mat. No.5 (1960), 165-179(Russian).

\bibitem{TWW}J. Traub, G. Wasilkowski, H. Wo\'{z}niakowski, Information Based Complexity, Academic Press, New York, 1988.

\bibitem{Vo}S.M. Voronin, N. Temirgaliev, On some applications of the Banach measure, Izv. Akad. Nauk. Kazakh. SSR Ser.Mat.Fiz. (1984) 8-11(Russian).

\bibitem{W1}H. Wang, Probabilistic and average widths of Sobolev spaces on compact two-point homogeneous spaces equipped with a Gaussian measure, Constr. Approx. 39 (2014) 485-516.

\bibitem{W2}H. Wang, Probabilistic and average linear widths of weighted Sobolev spaces on the ball equipped with a Gaussian measure, J. Approx. Theory 241 (2019) 11-32.

\bibitem{WZ}H. Wang, X. Zhai, Best approximation of functions on the ball on the weighted Sobolev space equipped with a Gaussian measure, J. Approx. Theory 162 (2010) 1160-1177.

\bibitem{Wo}H. Wo\'{z}niakowski, Average case complexity of linear multivariate problems, Part 1: Theory, Part 2: Applications, J. Complexity 8(4) (1992) 337-372, 373-392.

\bibitem{Y}D. Ylvisaker, Designs on random fields, in: Srivastava (Ed.), A survey of Statistical Design
and Linear Models, North-Holland, Amsterdam, (1975) 593-607.
\end{thebibliography}
\end{document}